\documentclass{article}
\usepackage{amssymb,amsmath,amsthm,enumerate,amsfonts}
\usepackage{color}
\makeatletter

\@addtoreset{equation}{section}
\makeatother
\renewcommand{\leq}{\leqslant}
\renewcommand{\geq}{\geqslant}
\renewcommand{\div}{\operatorname{div}}
\newcommand{\trace}{\operatorname{tr}}
\newcommand{\as}{\operatorname{as}}
\newcommand{\dist}{\operatorname{dist}}
\newcommand{\Id}{\operatorname{Id}}
\providecommand{\R}{\mathbb{R}}

\providecommand{\N}{\mathbb{N}}

\providecommand{\V}{\mathcal{V}}

\newcommand\calG{\mathcal{G}}

\newcommand\cW{\mathcal{W}}

\renewcommand{\leq}{\leqslant}
\renewcommand{\geq}{\geqslant}

\renewcommand{\div}{\operatorname{div}}
\newcommand{\curl}{\operatorname{curl}}
\newtheorem{Theorem}{Theorem}

\newtheorem{Proposition}[Theorem]{Proposition}
\newtheorem{Lemma}[Theorem]{Lemma}
\newtheorem{Remark}[Theorem]{Remark}
\author{Olivier Glass\footnote{Ceremade (UMR CNRS no. 7534),
Universit\'e Paris-Dauphine, 
Place du Mar\'echal de Lattre de Tassigny,
75775 Paris
FRANCE
},
Franck Sueur\footnote{Laboratoire Jacques-Louis Lions,
Universit\'e Pierre et Marie Curie - Paris 6,
4 Place Jussieu,
75005 Paris
FRANCE
}, }

\date{\today}
\title{On the motion of  a rigid body in a two-dimensional irregular ideal flow}
\begin{document}
\maketitle
\begin{abstract}
We consider the motion of a rigid body immersed in an ideal flow occupying the plane, with bounded initial vorticity. 
In that case there exists a unique corresponding solution which is global in time, in the spirit of the famous work by Yudovich for the fluid alone. 
We prove that if the body's boundary is Gevrey then the body's trajectory is Gevrey. 
This extends the previous work \cite{ogfstt} to a case where the flow is irregular.
\end{abstract}
\section{Introduction}
In this paper, we consider the movement of a rigid body immersed in a perfect incompressible fluid in the plane. The solutions that we will consider are weak solutions in the spirit of the solutions of Yudovich \cite{Yudovich} concerning the fluid alone. The goal of this paper is to study time regularity issues for the solid and fluid flows associated to such solutions of the system. \par
\ \par
Let us be more specific on the problem under view. We consider the motion of a rigid body which occupies at time $t$ the domain $\mathcal{S} (t) \subset \R^{2}$. The motion of this solid is rigid (and as we will see, driven by the pressure force on its boundary), so that ${\mathcal S}(t)$ is obtained by a rigid movement (that is a translation and a rotation) from its initial position ${\mathcal S}_{0}$, which is supposed to be a closed, simply connected domain in the plane with smooth boundary. In the rest of the plane, that is in the open set 
\begin{equation*}
\mathcal{F}(t) := \R^2 \setminus {\mathcal S} (t),
\end{equation*}
evolves a planar ideal fluid driven by the Euler equations. We denote correspondingly ${\mathcal F}_{0}:=\R^{2} \setminus {\mathcal S}_{0}$ the initial fluid domain. \par
The complete system driving the dynamics reads
\begin{gather}
\label{Euler1}
\frac{\partial u}{\partial t}+(u\cdot\nabla)u + \nabla p =0 \ \text{ for } \ x \in \mathcal{F}(t), \\
\label{Euler2}
\div u = 0 \ \text{ for } \  x \in \mathcal{F}(t) ,  \\
\label{Euler3}
u\cdot n =  u_\mathcal{S} \cdot n \ \text{ for } \ x\in \partial \mathcal{S}  (t),  \\
\label{Solide1} 
m  h'' (t) =  \int_{ \partial \mathcal{S} (t)} p n \, ds ,  \\
\label{Solide2} 
\mathcal{J}  \theta'' (t) =    \int_{ \partial   \mathcal{S} (t)} p  (x-  h (t) )^\perp  \cdot n \, ds , \\
\label{Eulerci2}
u |_{t= 0} = u_0 , \\
\label{Solideci}
h (0)= h_0 , \ h' (0)= \ell_0,  \theta (0) =  0 ,\ r  (0)=  r _0. 
\end{gather}
Here  $u=(u_1,u_2)$ and  $p$  denote the velocity and pressure fields defined on ${\mathcal F}(t)$ for each $t$,  $m>0$ and $ \mathcal{J}>0$ denote respectively the mass and the inertia of the body while the fluid  is supposed to be  homogeneous of density $1$, in order to simplify the equations (and without loss of generality). \par
When $x=(x_1,x_2)$ the notation $x^\perp$ stands for  $x^\perp =( -x_2 , x_1 )$, $n$ denotes the unit outward normal on $\partial \mathcal{F}(t)$, $ds$ denotes the integration element on the boundary  $\partial \mathcal{S}(t)$ of the body.
In the equations  (\ref{Solide1})  and (\ref{Solide2}), $h (t)$ is the position of the center of mass of the body,
$$\ell (t) := h' (t),$$
is the velocity of the center of mass and
$$r  (t):= \theta'(t),$$
is the angular velocity of the body. Accordingly, the solid velocity is given by
\begin{equation} \label{vietendue}
u_{{\mathcal S}}(t,x) := \ell(t) + r(t) (x-h(t))^{\perp}.
\end{equation}
Since  $\mathcal{S} (t)$ is a rigid body its position is obtained thanks to the rotation matrix 
\begin{eqnarray}
Q (t):= 
\begin{bmatrix}
\cos  \theta (t) & - \sin \theta (t)
\\  \sin  \theta (t) & \cos  \theta (t)
\end{bmatrix}.
\end{eqnarray}
More precisely the position $ \Phi^{{\mathcal S}} (t,x) \in \mathcal{S} (t)$  at 
the time $t$ of the point fixed to the body with an initial position $x$ is 
\begin{equation*}
\Phi^{{\mathcal S}} (t,x) := h (t) + Q (t)(x- h_0),
\end{equation*}
so that
\begin{equation*}
{\mathcal S}(t) = \Phi^{{\mathcal S}} (t,{\mathcal S}_{0}).
\end{equation*}
In Eulerian description  the velocity $u_\mathcal{S} (t,x) $ of the body  $\mathcal{S} $ at time $t$ at the position $x$ is
\begin{equation*}
u_\mathcal{S} (t,x) = {\partial_t}  \Phi^{{\mathcal S}} (t,   {(\Phi^{{\mathcal S}})}^{-1} (t,x) ) .
\end{equation*}
 The flow $ \Phi^{{\mathcal S}} $ corresponding to the solid is a rigid movement, that can be considered as a function of $t \in \R$ with values in the special Euclidean group $SE(2)$ of rigid movement in the plane. \par
The equations (\ref{Euler1}) and (\ref{Euler2}) are the incompressible Euler equations,  the  condition (\ref{Euler3}) means that the boundary is impermeable, the equations  (\ref{Solide1})  and (\ref{Solide2})  are  the Newton's balance laws for linear and angular momenta: the fluid acts on the body through pressure forces.\par
\ \par
A key quantity in the analysis is the vorticity 
\begin{equation*}
\omega  := \text{curl }  u = \partial_1 u_2 -   \partial_2 u_1,
\end{equation*}
which satisfies the transport equation:
\begin{eqnarray}
\label{eqvorty}
\partial_t \omega  + (u \cdot \nabla) \omega  = 0 .
\end{eqnarray}
One has the following result concerning the Cauchy problem for the above system, the initial position of the solid being given. This result describes both Yudovich and classical solutions.
\begin{Theorem} \label{ThmYudo}
For any $u_0 \in C^{0}(\overline{\mathcal{F}_0};\R^{2})$, $(\ell_0,r_0) \in \R^2 \times \R$, such that:
\begin{equation} \label{CondCompatibilite}
\div u_0 =0 \text{ in } {\mathcal F}_0 \ \text{ and } \  u_0   \cdot  n = (\ell_0 + r_0 (x-h_{0})^{\perp})   \cdot  n \text{ on } \partial \mathcal{S}_0,
\end{equation}
\begin{equation} \label{TourbillonYudo}
w_0 := \curl u_0  \in L_c^{\infty}(\overline{\mathcal F}_0),
\end{equation}
\begin{equation*}
\lim_{|x| \rightarrow +\infty} u_{0}(x) =0,
\end{equation*}
there exists a unique solution $(\ell,r,u)$ of \eqref{Euler1}--\eqref{Solideci} in $C^1 (\R; \R^2 \times \R) \times L^{\infty}_{loc}(\R; \mathcal{LL}({\mathcal F}(t)))$ with $\partial_{t} u, \nabla p \in L^{\infty}_{loc}(\R;L^{q}({\mathcal F}(t)))$ for any $q \in (1,+\infty)$. \par
Moreover such a solution satisfies that  for all $t>0$, $w(t):=\curl(u(t)) \in L^\infty_c(\overline{{\mathcal F}(t)})$, 
 and  $\| w(t,\cdot)\|_{L^{q}({\mathcal F}(t))}$ (for any $q \in [1,\infty]$), $\int_{{\mathcal F}(t)} w(t,x) \, dx$ and the circulation
 \begin{equation}
 \label{circu}
 \gamma := \int_{\partial {\mathcal S}(t)} u \cdot \tau \, ds,
\end{equation}
are preserved over time. \par
If moreover  $u_0 \in C^{\lambda + 1,\nu}({\mathcal{F}_0};\R^{2})$ for $\lambda$ in $\N$ and $\nu \in (0,1)$, then  $u$ is $L^{\infty}_{loc}(\R; C^{\lambda+ 1,\nu} ({\mathcal F}(t)))$.
\end{Theorem}
The notation $\mathcal{LL}(\Omega)$ refers to the space of log-Lipschitz functions on $\Omega$, that is the set of functions $f \in L^{\infty}(\Omega)$ such that
\begin{equation} \label{DefLL}
\| f \|_{\mathcal{LL}(\Omega)} := \| f\|_{L^{\infty}(\Omega)} + \sup_{x\not = y} \frac{|f(x)-f(y)|}{|x-y|(1+ \ln^{-}|x-y|)} < +\infty.
\end{equation}
On the other hand the notation  $C^{\lambda,\nu}(\mathcal{F}_{0})$ denotes the H\"older space, endowed with the norm:
\begin{align*}
\| u   \|_{  C^{\lambda,\nu} ( \mathcal{F}_{0}  ) } := 
\sup_{ |\alpha| \leqslant \lambda}   \big(  \|   \partial^\alpha u  \|_{L^\infty (  \mathcal{F}_{0} ) }  
+ \sup_{ x \neq y \in   \mathcal{F}_{0}  } \frac{ |\partial^\alpha u (x) -  \partial^\alpha u (y)| }{  |x - y|^\nu }  \big) <  + \infty  .
\end{align*}
Above, we used the abuse of notation $L^{\infty}(\R;X({\mathcal F}(t)))$ (resp. $C^{0}(\R;X({\mathcal F}(t)))$) where $X$ is a functional space; by this we refer to functions defined for almost each $t$ as a function in the space $X({\mathcal F}(t))$, and which can be extended as a function in $L^{\infty}(\R;X(\R^{2}))$ (resp. $C^{0}(\R;X(\R^{2}))$). \par
Theorem \ref{ThmYudo} is proven in \cite{shrinking}. Let us also mention that the existence and uniqueness of finite energy and classical solutions to the problem \eqref{Euler1}--\eqref{Solideci} has been tackled by Ortega, Rosier and Takahashi in \cite{ort1}. \par
\ \par
Consider a solution $(\ell,r,u)$ given by Theorem \ref{ThmYudo}. The corresponding fluid velocity field $u$ is log-Lipschitz in the $x$-variable; consequently there exists a unique flow map $\Phi^{{\mathcal F}}$ continuous from  $\R  \times{\mathcal F}_0 $ to $\mathcal{F}(t) $ such that
\begin{equation}
\label{flow}
\Phi^{{\mathcal F}} (t,x)  =  x + \int^{t}_0 u (s, \Phi^{{\mathcal F}} (s,x) ) ds .
\end{equation}
Moreover there exists $c >0$ such that for any $t$, the vector field  $ \Phi(t,\cdot) $ lies in the  H\"older space 
\begin{equation*}
C^{0 , \exp(-c|t|  \| \omega_0  \|_{L^{\infty}  ({\mathcal F}_0   )})}  ({\mathcal F}_0   ), 
\end{equation*}
and an example due to Bahouri and Chemin \cite{BahouriChemin} shows that this estimate is optimal. \par
For $M \geq 1$, we will denote by $\calG^{M} ((-T,T);E) $ the  Gevrey  space of order $M$ of smooth functions  $f:  (-T,T)  \rightarrow E$ with values in a Banach space $E$ such that for any compact set $K \subset (-T,T)$ there exist $L , C >0$ such that  for all $s \in \mathbb{N}$ and for all $t \in K$,
\begin{equation}
\label{croissderive}
\|  \partial_{t}^{s} f (t)  \|  \leq C L^s (s!)^{M}  .
\end{equation}
Let us recall that for $M=1$ the  Gevrey  space $\calG^{M} ((-T,T);E) $ is the space $C^{\omega}  ((-T,T);E) $ of real-analytic functions. \par
\ \par
The main result of this paper is the following.
\begin{Theorem}
\label{start4}
Assume that the boundary $\partial {\mathcal S}_{0}$ is Gevrey of order $M$ and that the assumptions of Theorem \ref{ThmYudo} are satisfied. Associate the solid and fluid flows $\Phi^{{\mathcal S}}$ and $\Phi^{{\mathcal F}}$ to the solution $(\ell,r,u)$.
Then for each $T>0$,
\begin{equation}
\label{Geg}
(\Phi^{{\mathcal S}}, \Phi^{{\mathcal F}} ) \in {\mathcal G}^{M+2} ( (-T,T); SE(2) \times C^{0 , \exp(-c T  \| \omega_0  \|_{L^{\infty}  ({\mathcal F}_0   )})}  ({\mathcal F}_0   )) .
\end{equation}
\end{Theorem}
\begin{Remark}
Since we will proceed by induction the proof of Theorem \ref{start4} also yields that assuming  that the boundary $\partial {\mathcal S}_{0}$ is  only  $C^{k+1,\nu}$, with $k  \in \N$ and $\nu \in (0,1)$, then the flow $(\Phi^{{\mathcal S}}, \Phi^{{\mathcal F}} )$ are   $C^k$ from $(-T,T)$ to $ SE(2) \times C^{0 , \exp(-c T  \| \omega_0  \|_{L^{\infty}  ({\mathcal F}_0   )} )}  ({\mathcal F}_0   )$.
 \end{Remark}
\begin{Remark}
Considering the particular case $M=1$, we see that when the boundary is real-analytic, the flows $(\Phi^{{\mathcal S}}, \Phi^{{\mathcal F}} )$ belong to the Gevrey space ${\mathcal G}^{3}$. This extends the results obtained in  \cite{gamblin} for a fluid filling the whole plane and in \cite{moi} for a fluid bounded by fixed boundaries. 
Moreover, mixing the techniques of the present paper and the ones of \cite{ogfstt}, one could prove that for strong solutions of the system, that is when $u_0 \in C^{\lambda + 1,\nu}({\mathcal{F}_0};\R^{2})$ in Theorem \ref{ThmYudo}, and when the boundary $\partial \mathcal{S}_0$ is real-analytic, then the flows $(\Phi^{{\mathcal S}}, \Phi^{{\mathcal F}} )$ belong to $C^{\omega}( (-T,T); SE(2) \times C^{\lambda+1 ,\nu}  ({\mathcal F}_0   )) .$ 
Actually it is expected that this also holds true in three dimensions locally in time. This is the equivalent in the case of an exterior domain of what is proven in \cite{ogfstt} in a bounded one.
\end{Remark}
Another way to express the result obtained in Theorem \ref{start4} is that we prove that for any $k \in \N$, for any  $T>0$, any $\tau \in (0,T)$, there exists $L>0$ such that for any  $t \in [-\tau,\tau]$,
\begin{eqnarray}
\label{PPPM}
\| \partial^{k +1 }_t \Phi^{{\mathcal F}} (t,\cdot )  \|_{ C^{0 , \exp(- c  T  \| \omega_0  \|_{L^{\infty}  ({\mathcal F}_0   )})}  ({\mathcal F}_0   )  }  + \| \ell^{(k)}  (t) \| + | r^{(k)}  (t)| 
& \leqslant  & L^{k +1 } (k!)^{M+2} .
\end{eqnarray}
%

%
%
%
%
%
%
%
\section{Preliminaries}
\label{SecPM}
In this section, we introduce the basic material that we use in order to prove Theorem \ref{start4}. \par
\ \par
\subsection{Basic material}
We introduce the distance to $\partial {\mathcal S}_{0}$:
\begin{equation*}
\rho(x) := \dist(x,\partial \mathcal{S}_0 ).
\end{equation*}
Using the assumption on ${\mathcal S}_{0}$, we deduce that there exists $c_\rho > 0$ such that on some bounded neighborhood ${\mathcal W}_0$ of the boundary $\partial \mathcal{S}_0 $ in $\overline{{\mathcal F}_{0}}$, one has for all $s \in \mathbb{N}$, 
\begin{equation}
\label{bordGevrey}
\| \nabla^s \rho \| \leq c_\rho^s \, (s!)^M ,
\end{equation}
as a function (on $\cW$) with values in the set of symmetric $s$-linear forms. Without loss of generality, we may assume that
\begin{equation*}
c_{\rho} >1.
\end{equation*}
Since the motion of the body is rigid we have that, for all $t$, the distance to ${\mathcal S}(t)$ in 
\begin{equation*}
{\mathcal W} (t) := \Phi^{{\mathcal S}} (t, {\mathcal W}_{0}),
\end{equation*}
is given  by
\begin{equation} \label{RhoRho}
\rho  (t,x)  =  \rho_0  ( (\Phi^{{\mathcal S}})^{-1}(t,x)).
\end{equation}
It will useful to have in mind the following form of the H\"older inequality: for any integer $k$, for any $ \theta := (s, \alpha)$ in
\begin{equation*}
\mathcal{A}_{k} := \{ \theta \in  \N^* \times  \N^s / \ 2 \leqslant s  \leqslant k+1  \text{ and }
\alpha := ( \alpha_1,\ldots, \alpha_s )  \in \N^s / \ | \alpha |    = k+1 - s \},
\end{equation*}
where the notation $ | \alpha |$ stands for $| \alpha | := \alpha_1 + \ldots+  \alpha_s$,
and for any $p \geqslant 1$, 
\begin{eqnarray}
\label{holder}
\left\| \prod_{i=1}^{s} f_i \right\|_{ L^{\frac{p}{k+1}} (\mathcal{F}(t)) }
\leqslant \prod_{i=1}^{s} \| f_i \|_{ L^{ \frac{p}{ \alpha_i +1} } (\mathcal{F}  (t) ) } .
\end{eqnarray}
\ \par
\noindent
{\bf Conventions.} We will use Einstein's repeated index convention. 
Given $A \in {\mathcal M}_{2}(\R)$, we denote by $\trace\{A\} $ its trace and by $\as\{A\}:=A-A^*$ its antisymmetric part.
Given $\varphi$ a smooth vector field, $\curl(\varphi)$ can be considered either as a scalar $\partial_1 \varphi_2 - \partial_2 \varphi_1$ or as a matrix, whose entry at the $k$-th row and $l$-th column is $[\curl(\varphi)]_{kl}= \partial_k \varphi_{l} - \partial_{l} \varphi_{k}$. The translation between the two is immediate. Also, $\nabla \varphi$ is the matrix $(\nabla \varphi)_{kl}=\partial_{k} \varphi_{l}$. Hence $\curl \varphi = \as \{ \nabla \varphi \}$. \par
Throughout this paper, we denote by $\N$ the set of nonnegative intergers and by $\N^{*}$ the set of positive integers. \par
\subsection{Added mass}
\label{Subsec:AD}
We will use the following decomposition of the pressure, which is the two-dimensional counterpart of \cite[Lemma 3]{ogfstt}.
\begin{Lemma}
Equation \eqref{Euler1} can  be written as
\begin{equation}
	Du=-\nabla \mu + \nabla \left( \Phi \cdot \begin{bmatrix}  \ell \\ r \end{bmatrix}' \right),
	\label{new2.0}
\end{equation}
with $\Phi:=(\Phi_{a})_{a=1,2,3}$, where the functions $\Phi_{a}=\Phi_{a}(t,x)$ and $\mu=\mu(t,x)$ are the solutions of the following problems: 
\begin{equation} \label{t1.3}
-\Delta {\Phi}_a = 0 \quad   \text{for}  \ x\in \mathcal{F}(t),
\end{equation}
\begin{equation} \label{t1.4}
{\Phi}_a(x) \rightarrow 0  \quad  \text{for}  \ x  \rightarrow \infty, 
\end{equation}
\begin{equation} \label{t1.5}
\frac{\partial {\Phi}_a}{\partial n}= K_a  \quad  \text{for}  \  x\in \partial \mathcal{S}(t),
\end{equation}
where
\begin{equation} \label{t1.6}
K_a := \left\{\begin{array}{ll} 
n_a & \text{if} \  a=1,2 ,\\ \relax
(x-h(t))^\perp \cdot n & \text{if} \ a=3,
\end{array}\right.
\end{equation}
and
\begin{equation} \label{t1.0}
-\Delta \mu = \trace\{ \nabla u \cdot \nabla u \} \ \text{ for } \ x\in \mathcal{F}(t), 
\end{equation}
\begin{equation} \label{t1.1}
\mu(x) \rightarrow 0 \ \text{ as }  \ x  \rightarrow \infty, 
\end{equation}
\begin{equation} \label{t1.2}
\frac{\partial \mu}{\partial n}= \sigma ,  \quad   \text{for}  \ x\in \partial \mathcal{S}(t), 
\end{equation}
where $u_\mathcal{S}=u_\mathcal{S}(t,x)$ is given by \eqref{vietendue} and where 
\begin{equation} \label{new0.1}
\sigma := \nabla^{2} \rho \, \{ u- u_\mathcal{S} , u- u_\mathcal{S} \} - n \cdot \big(r \left(2u-u_\mathcal{S}-\ell \right)^{\perp}\big).
\end{equation}
\end{Lemma}
Let us observe that 
\begin{equation*}
\Phi_a  (t,x)=   {\Phi}_{a} (0,  (\Phi^{\mathcal{S}})^{-1} (t,x)  ), 
\end{equation*}
so that the matrix
\begin{eqnarray*}
\mathcal{M}_2 = \begin{bmatrix} \int_{\mathcal{F} (t)} \nabla \Phi_a \cdot \nabla \Phi_b \ dx \end{bmatrix}_{a,b \in \{1,2,3\}}
\end{eqnarray*}
is time-independent.
It is also easy to see that the matrix $\mathcal{M}_{2}$ is symmetric and positive, as a Gram matrix (actually, when ${\mathcal S}_{0}$ is not a ball, it is even positive definite). Consequently the matrix
\begin{equation*}
\mathcal{M} :=\mathcal{M}_1+\mathcal{M}_2, \quad
 \mathcal{M}_1 := \begin{bmatrix} m \Id_2 & 0 \\ 0 & \mathcal{J}\end{bmatrix},
\end{equation*}
is symmetric and positive definite. Actually $\mathcal{M}$ is referred as the ``virtual  inertia tensor'', it incorporates the
 ``added  inertia tensor''  $\mathcal{M}_2$ which, loosely speaking,  measures how much the  surrounding fluid resists the acceleration as the body moves through it. Its relevance in that context is highlighted by the following property, which is the two-dimensional counterpart of \cite[Lemma 4]{ogfstt}.
\begin{Lemma}
The equations \eqref{Solide1}-\eqref{Solide2} can be written as
\begin{equation} \label{EvoMatrice}
\mathcal{M} \begin{bmatrix} \ell \\ r \end{bmatrix}' 
= \begin{bmatrix}  \displaystyle\int_{ \mathcal{F}(t)} \nabla \mu \cdot \nabla \Phi_a \, dx   \end{bmatrix}_{a \in \{1,2,3\}} .
\end{equation}
\end{Lemma}
%
%
%
%
%
%
%
%
%
\subsection{Regularity lemma}
We will use the following elliptic regularity estimate. It is a classical regularity estimate for the $\div$-$\curl$ elliptic system in $L^{p}$, except for what concerns the treatment of low frequencies, for which we rely on an approach due to T. Kato \cite{katoana}.
\begin{Lemma}
\label{triv}
There exists $c >0$ such that   for any $p \in (2,+\infty)$,
for any  smooth vector field $f \in L^{p}(\mathcal{F}_{0};\R^{2})$ satisfying
\begin{enumerate}[(i)]
\item $\div  f$ and $\curl f$ are in $L^{p}(\mathcal{F}_{0})$,
\item $\div  f = \partial_i \partial_h a_{ih} $ and $\curl f:= \partial_i \partial_h b_{ih} $ with the $a_{ih} $ and the $b_{ih} $ in $ L^{p}(\mathcal{F}_{0})$,
\item there exists $\phi$ in $W^{1,p}({\mathcal W})$ such that
$(n \cdot f ) |_{\partial \mathcal{F}_{0}} = \phi |_{\partial \mathcal{F}_{0}} $, 
\end{enumerate}
then $f \in W^{1,p} (\mathcal{F}_{0};\R^{2})$ and 
\begin{equation}
\label{reg1}
c \| f\|_{ W^{1,p} (\mathcal{F} ) }  \leq  p
\big( \| \div  f \|_{ L^{p}(\mathcal{F} ) } + \| \curl  f \|_{ L^{p}(\mathcal{F}_{0}) } + \|  \phi \|_{   W^{1,p}  ( \cW)  } \big) +
\| a_{ih} \|_{ L^{p}(\mathcal{F}_{0}) } + \| b_{ih} \|_{ L^{p}(\mathcal{F}_{0}) } 
+ \left| \int_{\partial \mathcal{S}_0 }  f \cdot \tau \, ds  \right| .
\end{equation}
\end{Lemma}
To obtain Lemma \ref{triv}, we are first going to prove the following lemma.
\begin{Lemma}
\label{triv0}
Let $\eta \in  C^{\infty} (\R^{2})$ such that $1-\eta \in C^{\infty}_{c}(\R^{2})$. 
There exists $c >0$ such that for any $p \in (2,+\infty)$,
for any smooth vector field $f \in L^{p}(\R^2;\R^{2})$ satisfying
\begin{enumerate}[(i)]
\item $\div f$ and $\curl f$ are in $ L^{p}(\R^2)$,
\item $\div  f = \eta \partial_i \partial_h a_{ih} $ and $\curl  f  := \eta \partial_i \partial_h b_{ih} $ with the $a_{ih} $ and the $b_{ih} $ in $ L^{p}(\R^2)$,
\end{enumerate}
then $f \in W^{1,p} (\R^2)$ and 
\begin{equation*}
c \| f\|_{ W^{1,p} (\R^2) }  \leq  p
 ( \| \div  f \|_{ L^{p}(\R^2) } + \| \curl  f \|_{ L^{p}(\R^2) } )  +  \| a_{ih} \|_{ L^{p}(\R^2) } + \| b_{ih} \|_{ L^{p}(\R^2) }    .
\end{equation*}
\end{Lemma}
\begin{proof}[Proof of Lemma \ref{triv0}]
According to the Biot-Savart formula, we have, for any $x \in \R^2 $,
\begin{eqnarray}
\label{for}
f(x) = \int_{\R^{2}} H(x-y)  \div  f (y) \, dy + \int_{\R^{2}} H^{\perp}(x-y)  \curl  f (y) \, dy ,
\end{eqnarray}
where 
\begin{equation*}
H(x) := \frac{x^{\perp}}{2\pi |  x|^{2} }. 
\end{equation*}
Here the difficulty lies in the estimate of the $L^{p}$ norm of $f$, since for the gradient we have the classical Calder\'on-Zygmund estimates:
\begin{eqnarray*}
c \| \nabla f\|_{ L^{p} (\R^2) }  \leq  p
 ( \| \div  f \|_{ L^{p}(\R^2) } + \| \curl  f \|_{ L^{p}(\R^2) } ) .
\end{eqnarray*}
To overcome this difficulty we will follow the strategy of the proof of \cite[Lemma 9.1]{katoana}.
Let us only deal with the first term in Eq. \eqref{for}; the second one can be tackled in a similar way.
We introduce  a smooth function  $\zeta \in C^{\infty}_{c}(\R^{2})$ such that $\zeta=1$  in the unit ball  $ B(0,1) $  and  $\zeta=0$ in the complementary set of $ B(0,2) $. 
Now we decompose the integral in two parts:
\begin{eqnarray}
\label{for2}
 \int_{\R^{2}} H(x-y)  \div  f (y) dy = \int_{\R^{2}} \zeta (x-y) H(x-y)  \div  f (y) dy + \int_{\R^{2}} (1- \zeta (x-y)) H(x-y)  \eta(y) \partial_i \partial_h a_{ih} (y) dy .
\end{eqnarray}
Using the classical Calder\'on-Zygmund theory we obtain that the norm in $W^{1,p} (\R^2)$ of the  first term in Eq. \eqref{for2} is bounded by
%
$c p \| \div  f \|_{ L^{p}(\R^2) } $, 
%
where the constant $c$ does not depend on $p >2$ (the dependence of the constants on $p$ is crucial here as well as in Yudovich's argument \cite{Yudovich}). Observe in particular that the kernel $\zeta H$ is integrable at infinity which yields the $L^p$ part of the previous estimate. \par
Now,  for the second part we integrate by parts twice, so that we get 
\begin{eqnarray*}
\int_{\R^{2}} (1- \zeta (x-y)) H(x-y)  \eta \partial_i \partial_h a_{ih} (y) dy  =  \int_{\R^{2}} \partial_i \partial_h \Big(  \eta  (y) (1- \zeta (x-y)) H(x-y)\Big) a_{ih} (y) dy .
\end{eqnarray*}
Since the kernel above is smooth we easily get by Young's inequality that the norm of this term in $W^{1,p} (\R^2)$  is bounded by  $C \| a_{ih} \|_{ L^{p}(\R^2) }  $.
This completes the proof of Lemma  \ref{triv0}.
\end{proof}
\begin{proof}[Proof of Lemma  \ref{triv}]
Let $R>0$ be large enough for $ \mathcal{S}_0 \subset B(0,R) $ and let us consider  a smooth function  $\eta$ such that $\eta=0$  in a neighborhood of $ \partial  \mathcal{S}_0$ and  $\eta=1$ in the complementary set of $ B(0,R) $. 
According to the previous lemma, the solution $F$ in  $W^{1,p} (\R^2)$  of 
\begin{eqnarray*}
\div F = \eta \div f \ \text{ and } \ \curl F = \eta \curl  f \ \text{ in } \ \R^2 ,
\end{eqnarray*}
satisfies the estimate \eqref{reg1} (Here we extend the functions $ \eta$, $f$, $a_{ih} $ and $b_{ih}$ by $0$ inside $ \mathcal{S}_0$ to be in position to apply Lemma  \ref{triv0}).
Now, the function
\begin{equation*}
\tilde{ f} := f-F,
\end{equation*}
has its divergence and rotational supported in the ball $B(0,R)$, and its circulation around ${\mathcal S}_{0}$ is given by
\begin{equation*}
\int_{\partial {\mathcal S}_{0}} \tilde{f} \cdot \tau \, ds = \int_{\partial {\mathcal S}_{0}} f \cdot \tau \, ds - \int_{{\mathcal S}_{0}} \div F \, dx.
\end{equation*}
Consequently there holds (extending $n$ in the neighborhood ${\mathcal W}_{0}$ of $\partial {\mathcal S}_0$):
\begin{eqnarray*}
\| \tilde{f} \|_{ L^{p} ({\mathcal F}_0) } & \leq & 
C p \Big( \| \div  \tilde{ f} \|_{ L^{p}({\mathcal F}_0) } + \| \curl  \tilde{ f} \|_{ L^{p}({\mathcal F}_0) }
+ \| \tilde{f} \cdot n \|_{ W^{1,p}  ({\mathcal W}_{0}) } \Big)
+ C \left| \int_{\partial \mathcal{S}_0 }  \tilde{f} \cdot \tau \, ds  \right| \\
& \leq &
C p \Big( \| \div  { f} \|_{ L^{p}({\mathcal F}_0) } + \| \curl { f} \|_{ L^{p}({\mathcal F}_0) }
+ \| \phi \|_{ W^{1,p}  ({\mathcal W}) } + \| a_{ih} \|_{ L^{p}(\mathcal{F}_{0}) } + \| b_{ih} \|_{ L^{p}(\mathcal{F}_{0}) }\Big)
+ C \left| \int_{\partial \mathcal{S}_0 }  f \cdot \tau \, ds  \right| .
\end{eqnarray*}
This concludes the proof of Lemma  \ref{triv}.
\end{proof}
\begin{Remark}
\label{AutreRegdivcurl}
Thanks to the invariance properties of the divergence and of the curl with respect to rotation and translation
the previous lemma holds for the domain $\mathcal{F} (t)$ with the same constant for any time $t$.
\end{Remark}
\subsection{Formal identities}
We will use some formal identities, which have already been obtained in \cite{ogfstt}, as a combinatorial refinement of the ones obtained by  Kato in \cite{katoana}. They concern the iterated material derivatives $(D^k u)_{ k \in \N^* }$, where 
\begin{equation*}
D := \partial_{t} + (u \cdot \nabla),
\end{equation*}
under the assumption that $(\ell,r,u)$ is a smooth solution of the above system. \par
We use the following notations: for $\alpha := ( \alpha_1,\ldots, \alpha_s )  \in \N^s$ we will denote $\alpha ! :=  \alpha_1 ! \ldots \alpha_s !$.
We denote for any integer $k$, 
\begin{equation*}
\mathcal{A}_{k} := \{ \theta:= (s, \alpha) \in  \N^* \times  \N^s / \ 2 \leqslant s  \leqslant k+1  \text{ and }
\alpha := ( \alpha_1,\ldots, \alpha_s )  \in \N^s / \ | \alpha |    = k+1 - s \},
\end{equation*}
where the notation $ | \alpha |$ stands for $| \alpha | := \alpha_1 + \ldots+  \alpha_s$. \par
\ \par
\noindent
Let us be given a smooth vector field $\psi$. 
We first recall some formal identities for $\div D^k \psi$, for $\curl D^k \psi$ of the iterated material derivatives $ (D^k \psi)_{ k \in \N^* }$ as combinations of the functionals 
\begin{gather} \label{DefsFetHNew}
f( \theta)  [u,\psi] : = \nabla D^{\alpha_1} u \cdot \ldots \cdot \nabla D^{\alpha_{s-1}} u\cdot \nabla D^{\alpha_s} \psi, 
\end{gather}
with $ \theta := (s, \alpha) \in \mathcal{A}_{k}$.
The precise statement is the following (see \cite[Prop. 6]{ogfstt}).
\begin{Proposition}\label{P1New}
For $k \in \N^*$, we have in $ \mathcal{F}(t) $
\begin{eqnarray}
\label{P1fNew}
\div D^k \psi=D^k\left(\div \psi\right) + \trace\left\{F^k [u,\psi]\right\} \text{ where } 
F^k [u,\psi] := \sum_{\theta   \in \mathcal{A}_{k}  } c^1_k (\theta  )  \, f(\theta)  [u,\psi], \\ 
\label{P2fNew}
\curl D^k \psi=D^k\left(\curl \psi\right) + \as\left\{G^k [u,\psi]\right\} \text{ where } 
G^k [u,\psi] :=  \sum_{\theta  \in \mathcal{A}_{k}  }  c^2_k (\theta ) \,  f(\theta)  [u,\psi],
\end{eqnarray}
and where for $i=1$, $2$ and $\theta = (s, \alpha)$, the $c^i_k (\theta )$ are integers satisfying 
\begin{equation} \label{Ci:3}
|c^i_k ( \theta) | \leqslant \frac{k!}{\alpha!}.
\end{equation}
\end{Proposition}
\begin{Remark}
In particular for $\psi=u$ and $k=1$ we obtain
\begin{equation*}
\div(Du) = - \div (\nabla \mu) = \trace(F^{1}[u,u]).
\end{equation*}
\end{Remark}
Now we recall some formal identities for normal traces on $\partial S(t)$ of the rigid body of iterated material derivatives $D^{k} \psi$, and for iterated material derivatives of the functions $K_i$ defined in \eqref{t1.6}. 
To a scalar $r \in \R$ we associate the matrix 
\begin{equation*}
\mathcal{R} (r) := r \begin{pmatrix} 0 & -1 \\ 1 & 0\end{pmatrix}. 
\end{equation*}
To any $\beta \in \N^s$ and $r \in C^{|\beta|}((-T,T);\R)$ we define the functional $\mathcal{R}_\beta [r]$ which associates to the time-dependent function $r$ the time-dependent rotation matrix
\begin{equation} \label{RComposes}
\mathcal{R}_\beta [r] := \mathcal{R} \left( r^{(\beta_1) }\right) \cdot \ldots \cdot  \mathcal{R} \left(r^{(\beta_s) }\right) .
\end{equation}
For any $s  \in  \N^*$,  we will use some multi-indices ${\bf s}' := (s'_1,...,s'_s )$ in $\N^s$. Then we will denote $s' := |{\bf s}'| =s'_1 + ... + s'_s $, $ ( \underline{\alpha}_1 ,..., \underline{\alpha}_{s} )$ will be in $\N^{s'_1} \times ... \times \N^{s'_s} $ and $\alpha := ( \underline{\alpha}_1 ,..., \underline{\alpha}_{s} , \alpha_{s' +1} , ...,\alpha_{s' +s } ) $ will be an element of $ \N^{ s' + s}$.
The  bricks of  the following formal identities will be the functionals, defined for smooth vector fields $\varphi$ and $\psi$ and a multi-index $\zeta :=(s,{\bf s}',\alpha)  \in  \N^* \times \N^s \times \N^{s+s'}$:
\begin{eqnarray}
\label{defh}
h (\zeta) [r,\varphi,\psi] :=   \nabla^{s} \rho  (t,x) \{ \mathcal{R}_{\underline{\alpha}_1 } [r] D^{\alpha_{s' +1} } \varphi ,...,   \mathcal{R}_{\underline{\alpha}_{s-1} } [r] D^{\alpha_{s' +s-1} } \varphi,\mathcal{R}_{\underline{\alpha}_s } [r] D^{\alpha_{s' +s} } \psi \}  .
\end{eqnarray}
In  (\ref{defh}) the term  $\mathcal{R}_{\underline{\alpha}_i } [r] $ should be omitted when  $s'_{i} := 0 $. 
We introduce the following set 
\begin{equation} \label{Eq:Defck}
\mathcal{B}_{k} := \{ \zeta =(s,{\bf s}',\alpha)  /  \ 2  \leqslant s + s' \leqslant k + 1 \text{ and } |\alpha| + s +s' =k+1\}.
\end{equation}
We have the following formal identity. 
\begin{Proposition}
\label{Prop:BodyDirichlet}
Given a smooth vector field $\psi$, for $k \in \N^*$, there holds on  the boundary $ \partial {\mathcal S}(t)$
\begin{eqnarray}
\label{P4fbody}
n\cdot D^k \psi = D^k\left(n\cdot \psi\right) +H^k [r,u-u_\mathcal{S} ,\psi] \text{ where } 
H^k [r,u-u_\mathcal{S} ,\psi] := \sum_{ \zeta  \in  \mathcal{B}_{k} } \ d^{1}_k (\zeta ) \  h (\zeta) [r,u-u_\mathcal{S} ,\psi]  , \\ 
\label{t7.4}
D^{k} K_i = \widetilde{H}^k [r,u-u_\mathcal{S} , \sigma_i] \text{ where }  \widetilde{H}^{k} [r,u-u_\mathcal{S} , \sigma_i] := \sum_{ \zeta  \in  {\mathcal{B}}_{k} } \ {d}^{2}_{k} (\zeta ) \  h (\zeta) [r,u-u_\mathcal{S} , \sigma_i] ,
\end{eqnarray}
where the $K_i$ are defined in \eqref{t1.6},
\begin{equation} \label{DefSigmai}
\sigma_i:= e_i  \text{ if } i=1,2, \text{ and }
\sigma_i:=  (x- h(t))^{\perp} \text{ if } i=3,  
\end{equation}
and where the $d^{j}_k (\zeta)$, $j=1,2$, are integers satisfying, for any $\zeta := (s,{\bf s}',\alpha)\in \mathcal{B}_k$,
\begin{equation} \label{Di}
|d^{j}_k (\zeta) | \leqslant \frac{ 3^{s+s'} k ! }{\alpha ! (s-1)!},
\end{equation}
\end{Proposition}
The proof of this proposition is completely identical to the proof of \cite[Prop. 8]{ogfstt} and is therefore omitted. \par
\ \par
We can also establish identities for the gradient $\nabla D^k\psi$ for a smooth scalar-valued function $\psi$:
\begin{Proposition}\label{P3New}
For $k\geq 1$, we have in the domain $\mathcal{F}(t)$
\begin{equation} \label{t3.4New}
D^{k} \nabla \psi = \nabla D^k\psi + K^{k}[u,\psi] , 
\end{equation}
where for $k\geq 1$,
\begin{equation} \label{RelKk}
K^{k}[u,\psi]:=-\sum_{r=1}^{k}  \dbinom{k}{ r} \nabla D^{ r-1} u \cdot D^{k-r} \nabla \psi.
\end{equation}
\end{Proposition}
The proof of this proposition is completely identical to the proof of \cite[Prop. 3.5]{katoana}, and is therefore omitted. \par
\subsection{Further formal identities} 
In this subsection, we give some other formal identities aimed at dealing with the far-field/low frequencies. This is inspired by \cite[Section 6]{katoana}.  Again, we assume here that $(\ell,r,u)$ is a smooth solution of the system. \par 
\ \par
We first recall the following commutation rules which allow to exchange $D$ and other differentiations. They are valid for $\psi$ a smooth scalar/vector field defined in the fluid domain:
\begin{eqnarray}
	D(\psi_1\psi_2)=(D\psi_1)\psi_2+\psi_1(D\psi_2) \label{t3.0},\\
	\partial_{k} (D\psi) -D(\partial_{k} \psi) = (\partial_{k} u_{j}) (\partial_{j} \psi), \label{t3.1}\\
	\div D\psi - D\div \psi = \trace \left\{(\nabla u)\cdot(\nabla \psi) \right\}, \label{t3.2}\\
	\curl D\psi - D\curl \psi = \as \left\{(\nabla u)\cdot(\nabla \psi) \right\}. \label{t3.3}
\end{eqnarray}
\subsubsection{An identity concerning the divergence}
The first formal identity of this section is given in the following statement.
\begin{Proposition}\label{further}
Let us consider $\psi$ a smooth vector field, and suppose that for some family $(\hat{\phi}_{ij})_{i,j=1\dots d}$ of smooth functions, one has in $ \mathcal{F}(t) $
\begin{equation}
\label{HypPhi}
\div(\psi) = \partial_{i} \overline{\phi}_{i} \ \text{ where } \   \overline{\phi}_{i}= \partial_{j} \hat{\phi}_{ij}.
\end{equation}
Then for $n \in \N^*$, we have in $ \mathcal{F}(t) $
\begin{equation*}
\div D^n \psi= \partial_i \partial_j \hat{\phi}_{ij}^n  [u,\psi],
\end{equation*}
where
\begin{equation} \label{FormuleHatPhi}
\hat{\phi}_{ij}^n [u,\psi] := \sum_{\xi   \in \hat{\mathcal{A}}^{1}_{n}  } 
\hat{c}^{n,1}_{i,j} (\xi  )  \,  \hat{g} (\xi)  [u,\psi] 
+ \sum_{\xi   \in \hat{\mathcal{A}}^{2}_{n}  }  \hat{c}^{n,2}_{i,j} (\xi  )  \, \check{g} (\xi) [u,\hat{\phi} ],
\end{equation}
with, for $\delta= 1$ or $2$,
\begin{gather}
\label{MatcalA}
{\mathcal{A}}^{\delta}_{n} :=  \{ \theta:= (s, \alpha) \in  \N^* \times  \N^s / \ 3-\delta \leqslant s  \leqslant n+1  \text{ and }
\alpha := ( \alpha_1,\ldots, \alpha_s )  \in \N^s / \ | \alpha |    = n+1 - s \}, \\
\label{8juin}
\hat{\mathcal{A}}_{n}^{\delta} :=   \{ \xi := (s,\alpha , k ,\lambda )   / \  (s,\alpha ) \in \mathcal{A}^{\delta}_{n}  \ \text{ and } (k ,\lambda)   \in \{ 1,...,d \}^{[s+\delta -3] +s}  \} , \\
\hat{g} (\xi)  [u,\psi] :=   D^{\alpha_1 } u_{\lambda_1} \cdot  \partial_{k_2} D^{\alpha_2 } u_{\lambda_2}  \cdot ...  \cdot  \partial_{k_{s-1 } } D^{\alpha_{s-1 }  } u_{\lambda_{s-1 } } \cdot D^{\alpha_{s }  } \psi_{\lambda_{s } } , \\
\check{g} (\xi)  [u,\hat{\phi}] :=   D^{\alpha_1 } u_{\lambda_1} \cdot  \partial_{k_2} D^{\alpha_2 } u_{\lambda_2}  \cdot ...  \cdot  \partial_{k_{s-1 } } D^{\alpha_{s-1 }  } u_{\lambda_{s-1 } } \cdot \partial_{k_{s}} D^{\alpha_{s }  } \hat{\phi}_{\lambda_{s}k_{s} } , 
\end{gather}
and where, for $\delta=1$, $2$, the $ \hat{c}^{n,\delta}_{i,j}(\xi)$ are integers satisfying 
\begin{equation} \label{anelka}
|   \hat{c}^{n,\delta}_{i,j}(\xi) | \leqslant 4^s \frac{n ! }{\alpha ! } .
\end{equation}
\end{Proposition}
We first state and prove two lemmas before establishing Proposition \ref{further}.
\begin{Lemma}
\label{LemmaFurtherDiv}
Under the hypothesis of Proposition \ref{further}, we have 
that for $n \in \N$,  for $x$ in $ \mathcal{F}(t) $
\begin{eqnarray}
\label{further1step1}
\div D^n \psi &=&  \partial_i  \overline{\phi}_i^{n} \\  
\label{further1step1bis}
\overline{\phi}_i^{n} &=&  \partial_j  \hat{\phi}^n_{i,j} ,
\end{eqnarray}
where the sequences $(\overline{\phi}_i^{n} )_{n \in \N}$ and $(\hat{\phi}^n_{i,j})_{n \in \N}$, are respectively defined by
\begin{eqnarray}
\label{further1step1bar}
\overline{\phi}_i^{0} =  \partial_{j} \hat{\phi}_{ij} \text{ and }\ 
\overline{\phi}_i^{n+1} &:=& D  \overline{\phi}_i^{n} - (\partial_{k} u^i ) (\overline{\phi}_k^{n} - D^n {\psi}_k ), \\
\label{further1step1chap}
\hat{\phi}^0_{i,j} =  \hat{\phi}_{ij}  \text{ and }\ 
\hat{\phi}_{i,j}^{n+1} &:=& D  \hat{\phi}^n_{i,j} -  \hat{\phi}^n_{i,k}  \cdot  \partial_{k} u_j + u^i (\partial_{k}  \hat{\phi}^n_{j,k} -D^n \psi_j ). 
\end{eqnarray}
\end{Lemma}
\begin{proof}
This is an induction argument. \par
\ \par
\noindent
$\bullet$ Let us first prove by iteration that \eqref{further1step1} holds true when the sequence $(\overline{\phi}_i^{n} )_{n \in \N}$ is defined by   \eqref{further1step1bar}. The case $n=0$ is precisely the hypothesis \eqref{HypPhi}. Let us now assume that \eqref{further1step1} and \eqref{further1step1bar} hold true for some $n \in \N$ and let us show the same for the rank $n+1$. \par
We first use the commutation rule \eqref{t3.2} to exchange $D$ and $\div$ to get 
\begin{eqnarray*}
\div D^{n+1} \psi &=& D \div D^n \psi + \trace (\nabla u \cdot \nabla D^n \psi ).
\end{eqnarray*}
Then we use \eqref{further1step1} and the commutation rule \eqref{t3.1} so that 
\begin{eqnarray*}
\div D^{n+1} \psi  &=& D \partial_{i}  \overline{\phi}_i^{n} + ( \partial_{i} u_k )(\partial_{k} D^n \psi_i ), \\
&=& \partial_{i} D \overline{\phi}_i^{n} -( \partial_{i} u_k )(\partial_{k} \overline{\phi}_i^{n} -  \partial_{k} D^n \psi_i ) ,
\end{eqnarray*}
Now, since $ \div u = 0$ we get
\begin{eqnarray*}
\div D^{n+1} \psi
&=&  \partial_{i} D \overline{\phi}_i^{n} - \partial_{k} ( ( \partial_{i} u^k ) ( \overline{\phi}_i^{n} - D^n  \psi_i )) , \\
&=&  \partial_{i} D \overline{\phi}_i^{n} - \partial_{i} ( ( \partial_{k} u^i ) ( \overline{\phi}_k^{n} - D^n  \psi_k )) ,
\end{eqnarray*}
after exchanging the dummy indices $i$ and $k$. This yields  \eqref{further1step1bar}. \par
\ \par
\noindent
$\bullet$ We go on with  \eqref{further1step1bis} and \eqref{further1step1chap}.  
We also proceed by iteration. The case $n=0$ is satisfied by assumption, see again  the hypothesis \eqref{HypPhi}. 
Now let us assume that \eqref{further1step1bis} holds true for some $n \in \N$  when the sequence  $(\hat{\phi}^{n}_{i,j} )_{n \in \N}$ is defined by    \eqref{further1step1chap}, and show the same at rank $n+1$. \par
Plugging  \eqref{further1step1bis} (at rank $n$) in \eqref{further1step1bar} we have
\begin{eqnarray*}
\overline{\phi}_i^{n+1} &=& D \partial_j \hat{\phi}^n_{i,j} - ( \partial_{k} u^i )  (\partial_{j} \hat{\phi}^n_{k,j} - D^n \psi_k ),
\end{eqnarray*}
Then we use the commutation rule \eqref{t3.1} and Leibniz' rule to get 
\begin{eqnarray*}
\overline{\phi}_i^{n+1} &=& \partial_j  D \hat{\phi}^n_{i,j} - (\partial_{j} u^k ) ( \partial_{k}  \hat{\phi}^n_{i,j} ) 
- \partial_{k} (u^i (\partial_{j} \hat{\phi}^n_{k,j} - D^n \psi_k ))
+  u^i  (\partial_{j} \partial_{k} \hat{\phi}^n_{k,j} - \partial_{k}  D^n \psi_k ) .
\end{eqnarray*}
This last term vanishes according to \eqref{further1step1} so that, using $\div u = 0$  and exchanging some dummy indices, we get
\begin{eqnarray*}
\overline{\phi}_i^{n+1} &=& \partial_{j} D \hat{\phi}^n_{i,j} 
- \partial_{k} ( (\partial_{j} u^k ) (  \hat{\phi}^n_{i,j} + u^i  ( \partial_{j}  \hat{\phi}^n_{k,j} - D^n \psi_k ) ), \\ 
&=& \partial_{j} D \hat{\phi}^n_{i,j}
-  \partial_{j} ( (\partial_{k} u^j ) (  \hat{\phi}^n_{i,k} + u^i  ( \partial_{k}  \hat{\phi}^n_{j,k} - D^n \psi_j ) ) , \\
&=&  \partial_{j} \hat{\phi}^{n+1}_{i,j} ,
\end{eqnarray*}
where $ \hat{\phi}^{n+1}_{i,j}$  is given by  \eqref{further1step1chap}.
\end{proof}
The second step of the proof of Proposition \ref{further} is the following lemma.
\begin{Lemma}
\label{LemmaFurtherDivComb}
Under the hypothesis of Proposition \ref{further}, we have 
that for $n \in \N$, we have in $ \mathcal{F}(t) $
\begin{eqnarray}\label{further1}
\overline{\phi}_i^{n} =  \sum_{\xi   \in \overline{\mathcal{A}}^{1}_{n}  } 
\overline{c}^{n,1}_{i,j} (\xi  )  \,  \overline{g} (\xi)  [u,\psi] 
+  \sum_{\xi   \in \overline{\mathcal{A}}^{2}_{n}  }  
\overline{c}^{n,2}_{i,j} (\xi  )  \, \tilde{g} (\xi) [u,\hat{\phi}]      , 
\end{eqnarray}
with
\begin{eqnarray} \label{DeCadix}
\overline{\mathcal{A}}_{n}^{\delta} :=   \{ \xi := (s,\alpha , k ,\lambda )   / \  (s,\alpha ) \in \mathcal{A}^{\delta}_{n}  \ \text{ and } (k ,\lambda)   \in \{ 1,...,d \}^{[s+ \delta-2]+s}  \} , \\ 
\label{Defg}
\overline{g} (\xi)  [u,\psi] :=  \partial_{k_{1 } }  D^{\alpha_1 } u_{\lambda_1} \cdot  \partial_{k_2} D^{\alpha_2 } u_{\lambda_2}  \cdot ...  \cdot  \partial_{k_{s-1 }} D^{\alpha_{s-1 }  } u_{\lambda_{s-1 } } \cdot D^{\alpha_{s }  } \psi_{\lambda_{s } } , \\
\label{DefTildeg}
\tilde{g} (\xi)  [u,\hat{\phi}] :=  \partial_{k_{1 } }  D^{\alpha_1 } u_{\lambda_1} \cdot  \partial_{k_2} D^{\alpha_2 } u_{\lambda_2}  \cdot ...  \cdot  \partial_{k_{s-1 }} D^{\alpha_{s-1 }  } u_{\lambda_{s-1 } } \cdot \partial_{k_{s}} D^{\alpha_{s }  } \hat{\phi}_{\lambda_{s } k_{s}} ,
\end{eqnarray}
and where, for $\delta=1$, $2$, the $ \overline{c}^{n,\delta}_{i,j}(\xi)$ are integers satisfying 
\begin{equation}
\label{cndelatij}
| \overline{c}^{n,\delta}_{i,j} (\xi)| \leqslant 2^{s} \frac{n ! }{\alpha ! } .
\end{equation}
\end{Lemma}
\begin{proof}
Let us first explain why formula \eqref{further1} holds without yet estimating the coefficients $\overline{c}^{n,\delta}_{i,j}(\xi)$. When $n=0$, it is clear that one can put $\overline{\phi}_{i}^{0}=\partial_{j} \hat{\phi}_{ij}$ in the form \eqref{further1}. Next, assuming that $\overline{\phi}_{i}^{n}$ has this form, we apply \eqref{further1step1bar}. It is obvious that the terms $(\partial_{k} u^i ) (\overline{\phi}_k^{n} - D^n {\psi}_k )$ can enter the right hand side of 
\eqref{further1}; only the term $D \overline{\phi}_i^{n}$ is not trivial. But actually, that $D \overline{\phi}_i^{n}$ can be put of the form \eqref{further1} is a consequence of the induction assumption, Leibniz's formula and the commutation rule \eqref{t3.1}. \par
It remains to make the claim quantitative, that is to say, to prove that the coefficients $\overline{c}^{n,\delta}_{i,j}(\xi)$ that we recover satisfy \eqref{cndelatij}. For $n=0$ the claim is trivial. Let us now discuss the passage from $n$ to $n+1$. To obtain $\overline{\phi}_i^{n+1}$, we apply relation \eqref{further1step1bar} and consider the resulting terms coming from those of \eqref{further1}. We first consider a term coming from the first sum of \eqref{further1}, that is, corresponding to $\delta=1$. The modifications needed in the case $\delta=2$ will be explained below. \par
\ \par
Let $ \xi \in \overline{\mathcal{A}}^{1}_{n+1}$ with $s \geqslant 3$.   When computing $\overline{\phi}_i^{n+1}$ with the relation  \eqref{further1step1bar} the term involving $\overline{g} (\xi)  [u,\psi]$ comes:
\begin{enumerate}
\item either from $D \overline{\phi}_i^{n}$: the term appears after applying $D$ to a function $\overline{g}(\xi')[u,\psi]$, $\xi' \in \overline{\mathcal{A}}^{1}_{n}$, and using the commutation rule \eqref{t3.1}. The function function $\overline{g}(\xi')[u,\psi]$ has one less material derivative or a term having a factor of the type $ \partial_{k_{j}} u_i$  less.  There are at most 
$$ \sum_{j=1}^s \alpha_{j} 2^{s} \frac{n ! }{\alpha ! }, $$
contributions of the first kind (one has added a materiel derivative to some $\overline{g}(\xi') [u,\psi]$), and at most 
$$ (s-1) 2^{s-1} \frac{n ! }{\alpha ! }, $$
of the second kind  (one has added a $(\partial_{k_{j}} u_i)$ factor to some $\overline{g}(\xi') [u,\psi]$).
\item either from $ (\partial_{k} u^i ) \overline{\phi}_k^{n} $. Again there are at most $ (s-1) 2^{s-1} \frac{n ! }{\alpha ! } $ contributions.
\end{enumerate}
Since $(\sum_{j=1}^s \alpha_{j}) + s-1=n+1$, in total, we have
\begin{equation*}
| \overline{c}^{n+1,\delta}_{i,j}(\xi) | \leqslant 2^{s} \frac{(n+1)!}{\alpha!}.
\end{equation*}
\ \par
Concerning the particular case  $s=2$, we have to take into account the additional term $(\partial_{k} u^i ) D^n {\psi}_k$ in relation  \eqref{further1step1bar}. But since  $\delta=1$, there are no terms corresponding to $s=1$, that is, in this case, there are only contributions ``of the first kind'' as referred to above. It follows that here the total contributions can be estimated from above by
$$ n 2^{2} n ! +1 \leqslant 2^{2}(n+1)! .$$
The case $\delta=2$ is similar, but two modifications are in order:
\begin{enumerate}
\item here $s$ starts from $s=1$,
\item there are no particular additional contributions for $s=2$.
\end{enumerate}
\end{proof}
We can now get back to the proof of Proposition \ref{further}. 
\begin{proof}[Proof of Proposition \ref{further}]
Again the case $n=0$ is clear, and we discuss the passage from $n$ to $n+1$. The induction relation on which we will rely is the following, directly deduced from Lemma \ref{LemmaFurtherDiv}:
\begin{equation*}
 \hat{\phi}_{i,j}^{n+1} = D  \hat{\phi}^n_{i,j} -  \hat{\phi}^n_{i,k}  \cdot  \partial_{k} u_j + u^i ( \overline{\phi}_j^{n} -D^n \psi_j ) .
\end{equation*}
We begin with $\delta=1$. A term $\hat{g} (\xi)  [u,\psi]$ with $\xi \in \hat{\mathcal A}^{1}_{n+1}$ (with $s \geqslant 3$) has been obtain through several ways:
\begin{enumerate}
\item either from $D \hat{\phi}_{i,j}^{n}$:  either from a term whose one factor has one less material derivative or from a term having a factor of the type $ \partial_{k_{j}} u_i$ less. There are at most
$$ \sum_{j=1}^s \alpha_{j} 4^{s} \frac{n ! }{\alpha ! }, $$
contributions of the first kind, and at most
$$ (s-1) 4^{s-1} \frac{n ! }{\alpha ! }, $$
of the second kind,
\item either from $ \hat{\phi}^n_{i,k}  \cdot  \partial_{k} u_j $, which gives again at most $ (s-1) 4^{s-1} \frac{n ! }{\alpha ! } $ contributions,
\item either from $u^i  \overline{\phi}_j^{n}$, giving a contribution at most $2^{s-1} \frac{n !}{\alpha !}$, according to Lemma \ref{LemmaFurtherDivComb}.
\end{enumerate}
Summing these contributions, this gives the conclusion. \par
\ \par
For what concerns the case $s=2$, we have again an additional term $u^{i} D^{n} \psi_{j}$, but here there is no contribution coming from $ \hat{\phi}^n_{i,k}  \cdot  \partial_{k} u_j $ because there are no terms corresponding to $s=1$. The conclusion follows as previously. \par
\ \par
Again, the case $\delta=2$ is {\it mutatis mutandis} the same, with $s$ starting from $s=1$, but no particular additional term for $s=2$.
\end{proof}
\subsubsection{An identity concerning the curl}
The second formal identity of this section is given in the following statement.
\begin{Proposition}
\label{furtherCurl}
Let us consider $\psi$ a smooth gradient vector field.
Then for $n \in \N^*$, we have in $ \mathcal{F}(t) $
\begin{equation} 
\label{further2Curl}
(\curl D^n \psi )_{kl} = \partial_i \partial_j \hat{\Gamma}_{ij,kl}^n  [u,\psi] ,
\end{equation}
where 
\begin{equation} \label{further2Curlb}
\hat{\Gamma}_{ij,kl}^n  [u,\psi] := \sum_{\xi   \in \hat{\mathcal{A}}^{1}_{n}  }  \hat{d}^{n}_{ij,kl} (\xi  )  \,  \hat{g} (\xi)  [u,\psi],
\end{equation}
where $\hat{\mathcal{A}^{1}_{n}}$ is defined in \eqref{MatcalA} and the $ \hat{d}^{n}_{i,j}(\xi)$ are integers satisfying 
\begin{equation}
\label{anelka8juin}
| \hat{d}^{n}_{ij,kl} (\xi)| \leqslant 4^s \frac{n ! }{\alpha ! } .
\end{equation}
\end{Proposition}
\begin{Remark}
That there is only one sum in \eqref{further2Curlb} while there were two in \eqref{FormuleHatPhi} is due to the fact that we suppose $\psi$ to be a gradient field, so that \eqref{further2Curl} is trivial for $n=0$, which simplifies the analysis. This is sufficient to our purpose.
\end{Remark}
Again we will need two preliminary lemmas before proving Proposition \ref{furtherCurl}.
The following lemma is the counterpart of Lemma \ref{furtherCurl}.
\begin{Lemma}
\label{LemmaFurtherCurl}
Under the hypothesis of Proposition \ref{furtherCurl}, we have 
that for $n \in \N$,  for $x$ in $ \mathcal{F}(t) $
\begin{eqnarray}
\label{CurlDnPsi}
(\curl D^n \psi )_{kl} &=& \partial_i  \overline{\Gamma}_{i,kl}^n  \\
\label{further2step1bis}
\overline{\Gamma}_{i,kl}^n  &=&  \partial_j  \hat{\Gamma}_{ij,kl}^n ,
\end{eqnarray}
where the sequences $(\overline{\Gamma}_{i,kl}^n )_{n \in \N}$, $( \hat{\Gamma}_{ij,kl}^n)_{n \in \N}$ are respectively defined by
\begin{gather}
\label{IterGamma} 
\overline{\Gamma}_{i,kl}^0 = 0 \ \text{ and } \ \overline{\Gamma}_{i,kl}^{n+1}   := D  \overline{\Gamma}_{i,kl}^{n} - (\partial_{k} u_i ) \overline{\Gamma}_{i,kl}^{n} + (\partial_{k} u_i ) D^n  \psi_l -  (\partial_{l} u_i ) D^n  \psi_k , \\ 
\label{IterGamma2}    
\hat{\Gamma}_{ij,kl}^0 = 0  \ \text{ and }  \ \hat{\Gamma}_{ij,kl}^{n+1} :=  D   \hat{\Gamma}_{ij,kl}^{n} -  (\partial_{k} u_j) \hat{\Gamma}_{ih,kl}^{n} -  u_i   \overline{\Gamma}_{j,kl}^{n} - \delta_{j,k} u_i D^n  \psi_l  - \delta_{j,l} u_i D^n  \psi_k  .
\end{gather}
\end{Lemma}
\begin{proof}[Proof of Lemma \ref{LemmaFurtherCurl}]
Let us first prove by iteration that \eqref{CurlDnPsi} holds true when the sequences  $(  \overline{\Gamma}_{i,kl}^n )_{n \in \N}$ are defined by \eqref{IterGamma}. The case $n=0$ holds true since  $\psi$ is a  gradient vector field by hypothesis.  Let us now assume that \eqref{CurlDnPsi}-\eqref{IterGamma} holds true for some $n \in \N$ and let us show the same for the rank $n+1$. \par
We first use the commutation rule \eqref{t3.3} to exchange $D$ and $\curl$ to get 
\begin{eqnarray}
(\curl D^{n+1} \psi )_{kl} &=&  (  D \curl ( D^n \psi ) +  \as\left\{ \nabla u  \cdot  \nabla  D^n \psi   \right\})_{kl} .
\end{eqnarray}
Then we use \eqref{CurlDnPsi}, the commutation rule \eqref{t3.1} and that $\div u = 0$ to get
\begin{eqnarray*}
(\curl D^{n+1} \psi )_{kl}   &=&  D \partial_{i}  \overline{\Gamma}_{j,kl}^{n} 
+ ( \partial_{k} u_h ) ( \partial_{h}   D^n \psi_l ) - ( \partial_{l} u_h ) ( \partial_{h}   D^n \psi_k ) \\
&=&  \partial_{i}  D  \overline{\Gamma}_{j,kl}^{n} -   ( \partial_{i} u_h )( \partial_{h} \overline{\Gamma}_{i,kl}^{n} ) -\partial_{h} (  ( \partial_{k} u_h ) D^n \psi_l -  ( \partial_{l} u_h )   D^n \psi_k ) \\
&=& \partial_{i}  D  \overline{\Gamma}_{j,kl}^{n} - 
\partial_{h} (  ( \partial_{i} u_h )  \overline{\Gamma}_{i,kl}^{n}  - ( \partial_{k} u_h ) D^n \psi_l -  ( \partial_{l} u_h )   D^n \psi_k ) .
\end{eqnarray*}
After exchanging the dummy indices $i$ and $h$ this yields 
\begin{eqnarray*}
(\curl D^{n+1} \psi )_{kl}  
&=&  \partial_{i} ( D  \overline{\Gamma}_{i,kl}^{n} -  (\partial_{h} u_i )  \overline{\Gamma}_{h,kl}^{n}
+  (\partial_{k} u_i ) D^n \psi_l - ( \partial_{l} u_i ) D^n \psi_k ) \\
&=&  \partial_{i}  \overline{\Gamma}_{i,kl}^{n+1} ,
\end{eqnarray*}
where $ \overline{\Gamma}_{i,kl}^{n+1}$ is given by  \eqref{IterGamma}. Now,
\begin{eqnarray*}
\overline{\Gamma}_{i,kl}^{n+1} &= &
D  \overline{\Gamma}_{i,kl}^{n}  -  (\partial_{h} u_i )  \overline{\Gamma}_{h,kl}^{n}
+  (\partial_{k} u_i ) D^n \psi_l - ( \partial_{l} u_i ) D^n \psi_k \\ 
&=&   D  \partial_{j}   \hat{\Gamma}_{ij,kl}^{n+1} 
 - \partial_{h} ( u_i  \overline{\Gamma}_{h,kl}^{n} ) + u_i \partial_{h}  \overline{\Gamma}_{h,kl}^{n}
 + \partial_{k} ( u_i   D^n \psi_l ) - u_i  \partial_{k} (  D^n \psi_l ) - \partial_{l} ( u_i   D^n \psi_k )
 +  u_i \partial_{l}  D^n \psi_k .
\end{eqnarray*}
Summing the third term with the fifth term and the last one yields $0$ so that 
\begin{eqnarray*}
\overline{\Gamma}_{i,kl}^{n+1} &= &
\partial_{j}  D  \hat{\Gamma}_{ij,kl}^{n+1}  - \partial_{h} ( \partial_{j} u_h  \hat{\Gamma}_{ij,kl}^{n+1} )
- \partial_{h} ( u_i \overline{\Gamma}_{h,kl}^{n} ) -  \partial_{k} ( u_i D^n \psi_l )
- \partial_{l} (u_i D^n \psi_k ) \\
&=& \partial_{j}  \hat{\Gamma}_{ij,kl}^{n+1}
\end{eqnarray*}
where $ \hat{\Gamma}_{ij,kl}^{n+1} $  is given by  \eqref{IterGamma2}.
\end{proof}
\begin{Lemma}
\label{LemmaFurtherCurlComb}
Under the hypothesis of Proposition \ref{furtherCurl}, we have 
that for $n \in \N$, we have in $ \mathcal{F}(t) $
\begin{eqnarray}
\label{further2}
\overline{\Gamma}^{n}_{i,kl} = \sum_{\xi   \in \overline{\mathcal{A}}^{1}_{n}  } 
\overline{d}^{n}_{i,kl} (\xi  )  \,  \overline{g} (\xi)  [u,\psi],
\end{eqnarray}
where $\overline{\mathcal{A}}^{1}_{n}$ is defined in \eqref{DeCadix}, $\overline{g} (\xi)$ is defined in \eqref{Defg} and where the $ \overline{d}^{n}_{i,kl}$ are integers satisfying 
\begin{equation} \label{EstDikl}
| \overline{d}^{n}_{i,kl} | \leqslant 2^s \frac{n ! }{\alpha ! } .
\end{equation}
\end{Lemma}
\begin{proof}[Proof of Lemma \ref{LemmaFurtherCurlComb}]
This is exactly the same proof as for Lemma \ref{LemmaFurtherDivComb} (even, it is simpler since there is only one type of terms here), except for what concerns $s=2$. For $s=2$, we see from \eqref{IterGamma} that a term  $\overline{g} (\xi)  [u,\psi]$ with $ \xi \in \overline{\mathcal{A}}^{1}_{n+1}$ comes either by adding a material derivative to a factor in a term of $\overline{\Gamma}^{n}_{i,kl}$, or from the additional terms $(\partial_{k} u_i ) D^n  \psi_l$ and $ -  (\partial_{l} u_i ) D^n  \psi_k $. Hence the total contributions can be estimated from above by
$$ n 2^{2} n ! +2 =  2^{2}(n+1)! - 4 n! +2 \leqslant 2^{2}(n+1)! .$$
\end{proof}
\begin{proof}[Proof of Proposition \ref{furtherCurl}]
One can observe that, except for what concerns $s=2$, the induction relations \eqref{IterGamma}-\eqref{IterGamma2} are exactly the same as \eqref{further1step1bar}-\eqref{further1step1chap}. The only difference consists in the number of additional terms which can be $2$ in \eqref{IterGamma2} when $j=k=l$. But the particular argument for $s=2$ works again: we have (at most) two additional terms, but the other contributions come only from adding a material derivative to an existing term. The conclusion follows as previously.
\end{proof}
%
%
%
%
%
%
%
%
%
%
%
%
%
%
%
\section{Proof of Theorem \ref{start4}}
\label{proof}
We are now in position to prove Theorem \ref{start4}. \par
\ \par
We proceed by regularization. Consequently, we will work from now on a smooth solution of the equation, without changing the notation. Since the estimates that we are going to establish are uniform with respect to the regularization parameter, the general result follows. We refer to  \cite{ogfstt} for more details on this step. Note in particular that we can use the formal identities of Section \ref{SecPM} that were derived under the assumption that $(\ell,r,u)$ is smooth. \par
The main argument is to prove by induction an estimate on the $k$-th material derivative of the fluid and the body velocities. \par
\subsection{Iteration}
Let 
\begin{equation} \label{nota}
k_0 \in \N^*, \ \nu \in (0,1), \   p_1  \geq \frac{ 2 }{1-\nu }
\text{ and } p_2   \geq  p_1 \cdot (k_0 + 1) .
\end{equation}
We are going to prove recursively that for $L>0$ large enough, for any integer $k  \leqslant  k_0$, 
\begin{equation} \label{proofit}
\|D^k u\|_{ W^{1, \frac{p_2}{k+1}} (\mathcal{F} (t) ) }   + \|\ell^{(k)}\| + | r^{(k)} |
\leq  \V_k   ,
\end{equation}
with
\begin{equation*}
\V_k := \frac{p_2^k (k!)^M L^k}{(k+1)^{2}} \V^{k+1},
\end{equation*}
where
\begin{equation*}
\V := \|u\|_{ W^{1, p_2} (\mathcal{F}  (t) ) } +\|\ell\| + |r| .
\end{equation*}
The norm on vectors of $\R^{2}$ (here $\ell$ and its derivatives) is the usual Euclidean one.  We will also the notation $ \|\cdot  \|$ for the associated  matrix norm. \par
\ \par
The inequality \eqref{proofit} is true for $k=0$. 
Now let us assume that Eq. \eqref{proofit} is proved up to $  k -1 \leqslant  k_0 -1$.  \par
\ \par
We will first prove the following proposition, which, under the induction hypothesis, allows to estimate the next iterated material derivative of the pressure field, as it decomposed in Paragraph \ref{Subsec:AD}.
\begin{Proposition} \label{pesti}
The functions $\Phi_a$ ($a=1,2,3)$ and $\mu$ satisfy the following assertions.
\begin{itemize}
\item There exists a positive constant $C_0=C_0({\mathcal S}_{0})$ such that
\begin{equation} \label{pressurephi0}
\sum_{1\leq a \leq 3}  \|\nabla \Phi_a \|_{ W^{1, p_2} (\mathcal{F} (t)  ) }   \leq C_0.
\end{equation}
\item There exists $\gamma_{1}$ a positive decreasing function with $\displaystyle {\lim_{L \rightarrow + \infty}} \gamma_{1}(L)=0$ such that if for all $j \leq  k -1$, 
\begin{equation}	\label{pressurehyp1} 
\|D^ j u\|_{ W^{1, \frac{p_2}{j+1}} (\mathcal{F}  (t) ) }  +\|\ell^{( j)}\| + \|r^{( j)}\| \leq   \V_{ j} ,
\end{equation}
then for all $1\leq j \leq k$, 
\begin{equation} \label{pressureclphi}
\sum_{1\leq a \leq 3}  \|D^j  \nabla \Phi_a \|_{ W^{1, \frac{p_2}{j+1}} (\mathcal{F} (t)  ) }  \leq  \gamma_1 (L) \frac{\V_{j}}{\V}.
\end{equation}
\item There exists a positive constant $C_0=C_0 ({\mathcal S}_{0})$ such that
\begin{equation} \label{pressuremudepart}
 \|\nabla \mu \|_{ W^{1, \frac{p_2}{2}} (\mathcal{F}  (t) ) } \leq C_0 \, \V^{2}.
\end{equation}
\item There exists $\gamma_{2}$ a positive decreasing function with $\displaystyle {\lim_{L \rightarrow + \infty}} \gamma_{2}(L)=0$ such that if for all $j \leq  k -1$, 
\eqref{pressurehyp1} holds true then for all $1\leq j \leq  k -1$, 
\begin{equation}\label{pressureclmu}
\|D^ j \nabla \mu \|_{ W^{1, \frac{p_2}{j+2}} (\mathcal{F}  (t) ) }   \leq \gamma_{2} (L) \V \, \V_{j} .
\end{equation}
\end{itemize}
\end{Proposition}
Proposition \ref{pesti} is proven in Subsections \ref{Proofpesti} and \ref{Proofpesti2}. \par
\ \par
The second proposition allows to propagate the induction hypothesis on the solution $(\ell,r,u)$ itself. 
\begin{Proposition}\label{PropositionSolideFluide}
There exist a positive decreasing functions $\gamma_{3}$ with $\displaystyle {\lim_{L \rightarrow + \infty}} \gamma_{3}(L)=0$ such that  if for all $j \leq k-1$, \eqref{pressurehyp1} holds,
then
\begin{equation} \label{MvtSolideFluide}
\| \ell^{(k)} \| + | r^{(k)}| + \|D^{k} u\|_{ W^{1, \frac{p_2}{k+1}} (\mathcal{F}  (t) ) }   \leq  \V_{k} \,  \gamma_{3}(L).
\end{equation}
\end{Proposition}
The  proof of Proposition \ref{PropositionSolideFluide} is given in Subsection \ref{ProofPropositionSolideFluide}.  \par
\ \par
Once Proposition \ref{PropositionSolideFluide} established, the claim that \eqref{proofit} holds true  up to $  k  \leqslant  k_0 $ is a direct induction argument. \par
\subsection{Proof of Proposition \ref{pesti}, functions $\nabla \Phi_{a}$}
\label{Proofpesti}
Recall that the functions $\nabla \Phi_a$ ($a =1,2,3$) defined by \eqref{t1.3}--\eqref{t1.5} satisfy 
\begin{equation} \label{new0.4}
\div \nabla \Phi_a=0 \quad \text{in} \ \mathcal{F}(t), 
\quad \curl \nabla \Phi_a =0  \quad \text{in} \ \mathcal{F}(t), 
\quad n\cdot \nabla \Phi_a = K_a \quad \text{on} \ \partial \mathcal{S}(t),
\quad \Phi_a(x) \rightarrow 0  \quad \text{for}  \ x  \rightarrow \infty .
\end{equation}
Then by applying the regularity lemma (Lemma \ref{triv}), we obtain immediately \eqref{pressurephi0}. \par
\ \par
Let us now prove the second point of Proposition \ref{pesti} by induction on $j$.
Let us therefore assume that the estimate \eqref{pressureclphi} holds true up to rank $j-1$ and  prove that then it is also true at rank $j$.
By applying $D^{j}$ to \eqref{new0.4} and by using Propositions \ref{P1New},  \ref{Prop:BodyDirichlet}, \ref{P3New},  \ref{further} (with $\psi = \nabla \Phi_a $, which satisfy the assumption \eqref{HypPhi} with  identically vanishing functions $\hat{\phi}_{ij}$) and \ref{furtherCurl} (also with $\psi = \nabla \Phi_a $) we obtain that $D^{j} \nabla \Phi_a$ satisfies the following relations
\begin{gather}
\div D^{j} \nabla \Phi_a =\trace\left\{F^{j} [u,\nabla \Phi_a]\right\} \quad \text{in} \ \mathcal{F}(t),  \\
\label{new0.5}
\curl D^{j}\nabla \Phi_a =\as\left\{G^{j} [u,\nabla \Phi_a  ]\right\} \quad \text{in} \ \mathcal{F}(t), \\
\label{new0.3}
n\cdot D^{j}\nabla \Phi_a = D^{j} K_a +H^{j}[r,u-u_\mathcal{S} ,\nabla \Phi_a] \quad \text{on} \ \partial \mathcal{S}(t)  \\
\div D^j  \nabla \Phi_a  = \partial_i \partial_h a^j_{ih}  [u, \nabla \Phi_a ]  \quad \text{in} \ \mathcal{F}(t),  \\
\curl D^j  \nabla \Phi_a   = \partial_i \partial_h b^j_{ih}  [u, \nabla \Phi_a ]   \quad \text{in} \ \mathcal{F}(t) , \\
\label{circu8juin}
\int_{\partial \mathcal{S} (t) }   D^{j} \nabla  \Phi_a  \cdot \tau \, ds  =
\int_{\partial \mathcal{S} (t) }  K^{j} [u, \Phi_a ]   \cdot \tau \, ds ,
\end{gather}
where
\begin{gather*}
a^j_{ih} [u, \nabla \Phi_a ] := \sum_{\xi \in \hat{\mathcal{A}}^{1}_{j} }  a^j_{ih} (\xi)  \,  \hat{g} (\xi) [u, \nabla \Phi_a], \\ 
b^j_{ih} [u, \nabla \Phi_a ] := \sum_{\xi \in \hat{\mathcal{A}}^{1}_{j} }  b^j_{ih} (\xi)  \,  \hat{g} (\xi) [u, \nabla \Phi_a].
\end{gather*}
Above the set $\hat{\mathcal{A}}^{1}_{j} $ and the functionals  $ \hat{g} (\xi)$ are  the one defined in  \eqref{8juin}; 
and the integer coefficients $   a^j_{ih} (\xi  ) $ and $ b^j_{ih} (\xi  ) $ satisfy respectively  the estimates
\begin{eqnarray*}
\label{anelkaBIS}
| a^j_{ih} (\xi) | \leqslant 4^s \frac{j!}{\alpha!},
\quad
| b^j_{ih} (\xi) | \leqslant 4^s \frac{j!}{\alpha!}.
\end{eqnarray*}
\subsubsection{Estimate of $F^{j} [u,\nabla \Phi_a]$ and $G^{j} [u ,\nabla \Phi_a  ] $}
Applying the H\"older inequality \eqref{holder} to  $f(\theta)  [u ,\nabla \Phi_a]$ (whose definition is given in \eqref{DefsFetHNew}) for $\theta  \in \mathcal{A}_{j} $, we obtain that 
\begin{equation*}
 \| f( \theta)  [u,\nabla \Phi_a]  \|_{ L^{ \frac{p_2}{j+1}} (\mathcal{F}(t) ) }
 \leq  \prod_{i=1}^{s-1}  \| D^{\alpha_i } u \|_{ W^{1,\frac{p_2}{ \alpha_i +1}}   (\mathcal{F}(t) ) }    \| D^{\alpha_s }\nabla \Phi_a  \|_{ W^{1,\frac{p_2}{ \alpha_s +1}}   (\mathcal{F}(t) ) }   .
\end{equation*}
Using  the induction hypothesis and since for  $\theta  \in \mathcal{A}_{j} $, $| \alpha | = j+1 - s$, we have
\begin{equation*}
 \| f( \theta)  [u,\nabla \Phi_a]  \|_{ L^{ \frac{p_2}{j+1}} (\mathcal{F}(t) ) }
 \leq L^j  \,     \V^{j}            \, (\alpha !)^M
 L^{1-s}  p_2^{| \alpha | }   \prod_{i=1}^s  \frac{1}{ (1+\alpha_i )^2}.
\end{equation*}
Now thanks to Proposition \ref{P1New}, we obtain
\begin{eqnarray*}
\| F^j [u,\nabla \Phi_a] \|_{ L^{\frac{p_2}{j+1}}(\mathcal{F}(t)) }
\leq j! L^j \, \V^{j} \, \sum_{s=2}^{j+1} \ L^{1-s}  \, 
\sum_{\substack{{\alpha \text{ s.t.}} \\ {|\alpha| = j+1-s}}}  \,   (\alpha !)^{M-1} p_2^{ | \alpha |} \prod_{i=1}^s  \frac{1}{ (1+\alpha_i )^2} .
\end{eqnarray*}
When  $\theta  \in \mathcal{A}_{j} $, $2  \leq s \leq  j+1$ and $| \alpha | = j+1 - s$, then $| \alpha |  \leq j-1$ so that 
\begin{eqnarray}\label{new0.0}
\| F^j  [u,\nabla \Phi_a]  \|_{ L^{ \frac{p_2}{j+1}} (\mathcal{F}(t) ) }
 \leq  p_2^{j-1}  (j!)^M L^j  \,  \V^{j}  \,  \sum_{s=2}^{j+1} \  \frac{L^{1-s}}{j^{M-1}}  \, 
\sum_{ \alpha  /  \, | \alpha | = j+1 - s }  \,  \prod_{i=1}^s  \frac{1}{ (1+\alpha_i )^2}  .
\end{eqnarray}
We now use the following lemma (cf. \cite[Lemma 7.3.3]{cheminsmf}).  
%
%
\begin{Lemma} \label{LemmeCheminSMF}
For any couple of positive integers $(s,m)$ we have
\begin{equation}\label{DefUpsilon}
\sum_{\substack{{\alpha \in (\N^* )^{s}} \\ {|\alpha|=m} }}  \Upsilon(s,\alpha) \leq \frac{20^{s}}{(m+1)^{2}}, \text{ where }
\Upsilon(s,\alpha):=\prod_{i=1}^s \frac{1}{(1+\alpha_{i})^2}.
\end{equation}
\end{Lemma}
We deduce from \eqref{new0.0} and from the above lemma that
\begin{equation*}
\| F^j  [u,\nabla \Phi_a]  \|_{ L^{ \frac{p_2}{j+1}} (\mathcal{F}(t) ) }
\leq  p_2^{j-1} \frac{ (j!)^M L^j }{(j+1)^2} \V^{j}  \,  \sum_{s=2}^{j+1} \  \frac{L^{1-s}}{j^{M-1}} \,  20^{s}   \, \frac{(j+1)^2}{ (j-s+2 )^2}.
\end{equation*}
We obtain the same bound on $\| G^j [u,\nabla \Phi_a] \|_{ L^{ \frac{p_2}{j+1}} (\mathcal{F}(t)) }$ by using \eqref{P2fNew} instead of \eqref{P1fNew}.
\subsubsection{Estimate of $H^j [r,u-u_\mathcal{S},\nabla \Phi_a]$}
To estimate the body velocity $u_\mathcal{S} $ in $h (\zeta) [r,u-u_\mathcal{S} ,\nabla \Phi_a]$, we will use the following result, which is the two-dimensional counterpart of \cite[Lemma 8]{ogfstt} (with a different norm which has no importance here, since we estimate the solid velocity which belongs to a finite-dimensional space).
\begin{Lemma}\label{bodyvesti}
Under the same assumptions as Proposition \ref{pesti}, there exists a geometric constant $C > 1$ such that for any $m \leq k$
\begin{equation}
\|D^{m} \, u_\mathcal{S} \|_{ W^{1, \frac{p_2}{m+1}} ( \cW (t)  ) }  \leq   C   \V_{m}.
\end{equation}
\end{Lemma}
Applying \eqref{holder} to $h(\zeta)$ (whose formula is given in \eqref{defh}) and using \eqref{bordGevrey}, we obtain:
\begin{multline*}
\| h (\zeta) [r,u-u_\mathcal{S} ,\nabla \Phi_a]  \|_{  W^{1, \frac{p_2}{j+1}} (\cW  (t)  ) } \\
\quad \leq c_\rho^s   \,    (s!)^M  
\Big( \prod_{i=1}^s  \prod_{l=1}^{s'_i} | r^{(j)} |  \Big)
\Big(\prod_{i=1}^{s-1}   \|  D^{\alpha_{s' + i} } (u -u_\mathcal{S} ) \|_{ W^{1,\frac{p_2}{ \alpha_{s' + i} +1}}   (\cW (t) ) } \Big)
\|  D^{\alpha_{s' + s} }  \nabla \Phi_a \|_{ W^{1,\frac{p_2}{ \alpha_{s' + s} +1}}   (\mathcal{F}(t) ) } .
\end{multline*}
By using  the induction hypothesis and Lemma \ref{bodyvesti}, we have
\begin{equation*}
 \| h (\zeta) [r,u-u_\mathcal{S} ,\nabla \Phi_a]  \|_{  W^{1, \frac{p_2}{j+1}} (\cW  (t)  ) }
 \leq L^j  \, \V^{j}  \, (\alpha !)^M  (s!)^M
 L^{1-s} c_\rho^s   \,    p_2^{| \alpha | }   \prod_{i=1}^s  \frac{1}{ (1+\alpha_i )^2}.
\end{equation*}
Thanks to Proposition \ref{P1New} and Lemma \ref{LemmeCheminSMF} we obtain
\begin{eqnarray}
\label{new1.4}
\| H^j  [ r,u-u_\mathcal{S} ,\nabla \Phi_a] \|_{W^{1,\frac{p_2}{j+1}} (\cW  (t)  )} \leq  p_2^{j-1}
\frac{ (j!)^M L^j }{(j+1)^2} \V^{j}  \sum_{s=2}^{j+1} \ s^M  L^{1-s}  \, c_\rho^s   \,  20^{s} \,   \frac{(j+1)^2}{ (j-s+2 )^2}  .
\end{eqnarray}
\subsubsection{Estimate of $K^j  [u , \Phi_a]$}
Applying H\"older's inequality to the terms appearing in \eqref{RelKk} yields for any $j \geqslant 2$, 
\begin{eqnarray*}
\|  K^j  [u, \Phi_a]  \|_{ L^{ \frac{p_2}{j+1} } (\mathcal{F}(t) ) }  \leqslant 
\sum_{s=1}^{j-1} \  \dbinom{j-1}{s}  \| D^{s-1} u \|_{ W^{1,\frac{p_2}{s} }  (\mathcal{F}(t) )}
\| D^{j-s} \nabla  \Phi_a \|_{ W^{1,\frac{p_2}{j+1-s} } (\mathcal{F}(t) )}  .
\end{eqnarray*}
By using the induction hypothesis   we get 
\begin{eqnarray*}
\|  K^j  [u, \Phi_a]  \|_{ L^{ \frac{p_2}{j+1} } (\mathcal{F}(t) ) } 
\leqslant  p_2^{j-1} \frac{ (j!)^M L^j }{(j+1)^2}  \V^{j+1}  \,  L^{-1} 
\sum_{s=1}^{j-1} \  \left( \frac{j-s }{j s} \right)^M  \left( \frac{j+1 }{s(j-s+1)}\right)^2 ,
\end{eqnarray*}
and the same estimate holds true for 
$\int_{\partial \mathcal{S} (t) }   (D^{j} \nabla  \Phi_a)  \cdot \tau \, ds $ up to a multiplicative geometric constant,  thanks to \eqref{circu8juin}.
\subsubsection{Estimate of $a^j_{ih}  [u, \nabla \Phi_a ] $ and of $b^j_{ih}  [u, \nabla \Phi_a ] $}
Applying the H\"older inequality \eqref{holder} to the definition of  $ \hat{g} (\xi)  [u, \nabla \Phi_a ] $ in \eqref{8juin},  for $\xi  \in \hat{\mathcal{A}
}^{1}_{j}$,   yields that 
\begin{equation*}
 \|  \hat{g} (\xi)  [u,\nabla \Phi_a]  \|_{ L^{ \frac{p_2}{j+1}} (\mathcal{F}(t) ) }
 \leq  \prod_{i=1}^{s-1}  \| D^{\alpha_i } u \|_{ W^{1,\frac{p_2}{ \alpha_i +1}}   (\mathcal{F}(t) ) }    \| D^{\alpha_s }\nabla \Phi_a  \|_{ W^{1,\frac{p_2}{ \alpha_s +1}}   (\mathcal{F}(t) ) }   .
\end{equation*}
Then we proceed as for  $F^{j} [u,\nabla \Phi_a]$, and obtain
\begin{equation*}
\| a^j_{ih}  [u, \nabla \Phi_a ]  \|_{ L^{ \frac{p_2}{j+1}} (\mathcal{F}(t) ) }
\leq  p_2^{j-1} \frac{ (j!)^M L^j }{(j+1)^2} \V^{j}  \,  \sum_{s=2}^{j+1} \  \frac{L^{1-s}}{j^{M-1}} \,  80^{s}   \, \frac{(j+1)^2}{ (j-s+2 )^2}.
\end{equation*}
The analysis is the same for $b^j_{ih}  [u, \nabla \Phi_a ] $.
\subsubsection{Conclusion}
\label{ccl}
It remains to gather the above estimates. We apply Lemma  \ref{triv} (observing that, thanks to \eqref{nota}, we have  $ \frac{ p_2 }{ k+1}  > 2$) and we use the previous estimates  to get    \eqref{pressureclphi} at rank $j$, with
\begin{equation*}
\gamma_{1}(L) :=  C({\mathcal S}_{0}) \, \sup_{j \geq 1} \left( \sum_{s=2}^{j+1} \  \frac{L^{1-s}}{j^{M-1}} \,  80^{s} c_{\rho}^{s} \, \frac{(j+1)^2}{ (j-s+2 )^2} \right) .
\end{equation*}
\subsection{Proof of Proposition \ref{pesti}, function $\nabla \mu$}
\label{Proofpesti2}
We now turn to the claims concerning $\mu$. The function $\nabla \mu$ defined by \eqref{t1.0}--\eqref{t1.2} satisfies
\begin{gather*}
\div \nabla \mu=-\trace\left\{F^1[u,u]\right\}= -\trace\left\{\nabla u \cdot \nabla u \right\}, \quad \curl \nabla \mu=0 \quad \text{in} \ \partial \mathcal{F}(t), \\
\quad n\cdot \nabla \mu =  \sigma  \quad \text{on} \ \partial \mathcal{S}(t),
\end{gather*}
where $\sigma$ is defined by \eqref{new0.1}. Hence \eqref{pressuremudepart} follows again from Lemma \ref{AutreRegdivcurl}. \par
\ \par
By applying $D^{j}$ to \eqref{new0.4} and by using Propositions \ref{P1New},  \ref{Prop:BodyDirichlet}, \ref{P3New},  \ref{further} with $\psi = \nabla\mu $, which satisfies the assumption \eqref{HypPhi} with 
\begin{equation*}
\hat{\phi}_{ih} := -u_i u_h,
\end{equation*}
and Proposition \ref{furtherCurl} (also with $\psi = \nabla\mu $) we obtain that $D^{j} \nabla\mu$ satisfies the following relations
\begin{gather*}
\div D^{j} \nabla\mu =  - D^{j} \trace\left\{F^1[u,u]\right\} + \trace\left\{F^{j} [u,\nabla\mu]\right\} \quad \text{in} \ \mathcal{F}(t), \\
\label{new0.5mu}
\curl D^{j}\nabla\mu =\as\left\{G^{j} [u,\nabla\mu  ]\right\} \quad \text{in} \ \mathcal{F}(t), \\
\label{new0.3mu}
n\cdot D^{j}\nabla\mu = D^{j} \sigma +H^{j}[r,u-u_\mathcal{S} ,\nabla\mu] \quad \text{on} \ \partial \mathcal{S}(t)  \\
\div D^j  \nabla\mu  = \partial_i \partial_h a^j_{ih}  [u, \nabla\mu ] \quad \text{in} \ \mathcal{F}(t),  \\
\curl D^j  \nabla\mu   = \partial_i \partial_h b^j_{ih}  [u, \nabla\mu ]   \quad \text{in} \ \mathcal{F}(t) , \\
\label{circu8juinmu}
\int_{\partial \mathcal{S} (t) }   D^{j} \nabla \mu  \cdot \tau \, ds  =
\int_{\partial \mathcal{S} (t) }  K^{j} [u,\mu ]   \cdot \tau \, ds ,
\end{gather*}
where
\begin{gather*}
a^j_{ih}  [u, \nabla\mu ] :=
\sum_{\xi \in \hat{\mathcal{A}}^{1}_{j} } a^{j,1}_{ih} (\xi) \, \hat{g} (\xi) [u, \nabla\mu ]
+ \sum_{\xi \in \hat{\mathcal{A}}^{2}_{j} } a^{j,2}_{ih} (\xi) \, \hat{g} (\xi) [u, \hat{\phi} ],  \\
b^j_{ih}  [u, \nabla\mu ] := \sum_{\xi \in \hat{\mathcal{A}}^{1}_{j} } b^j_{ih} (\xi) \, \hat{g} (\xi) [u, \nabla\mu],
\end{gather*}

and the integer coefficients $a^{j,1}_{ih}(\xi)$, $a^{j,2}_{ih}(\xi)$ and $b^j_{ih}(\xi) $ satisfy  the estimate \eqref{anelkaBIS}. \par
The proof that the validity of \eqref{pressurehyp1} for $j \leq k-1$ implies the one of \eqref{pressureclmu} for $1\leq j\leq k-1$ is completely similar to the equivalent proof for $\Phi_{a}$. 

\subsection{Proof of Proposition  \ref{PropositionSolideFluide}}
\label{ProofPropositionSolideFluide}
We now turn to the proof of Proposition \ref{PropositionSolideFluide}. Under the same assumption that \eqref{pressurehyp1} is valid for all $j \leq k-1$, we first prove
\begin{equation} \label{MvtSolide}
\| \ell^{(k)} \| + | r^{(k)} | \leq  \V_{k} \,  \gamma_{4}(L),
\end{equation}
and then prove 
\begin{equation} \label{MvtFluide}
\| D^{k} u\| \leq   \,   \V_{k} \,   \gamma_{5}(L),
\end{equation}
for positive decreasing functions $\gamma_{4}, \gamma_{5}$ with $\displaystyle {\lim_{L \rightarrow + \infty}} \gamma_{4}(L)+ \gamma_{5}(L)=0$. \par
\ \par
In order to prove  \eqref{MvtSolide} it suffices to  differentiate the equations \eqref{EvoMatrice} $k$ times with respect to the time (recall that the matrix ${\mathcal M}$ is constant):
\begin{eqnarray*}
\mathcal{M}\begin{bmatrix} \ell \\ r \end{bmatrix}^{(k)}
&=&  \begin{bmatrix}  \displaystyle\int_{ \mathcal{F}(t)} \nabla \mu \cdot \nabla \Phi_a \, dx   \end{bmatrix}_{a \in \{1,2,3\}}^{(k-1)} ,
\\  &\leqslant& \begin{bmatrix}  \displaystyle\int_{ \mathcal{F}(t)} D^{k-1} (\nabla \mu \cdot \nabla \Phi_a ) \, dx   \end{bmatrix}_{a \in \{1,2,3\}} ,
\\  &\leqslant&  \sum_{i=0}^{k-1} \binom{k-1}{i}  
\begin{bmatrix}  \displaystyle\int_{ \mathcal{F}(t)}  (D^{i} \nabla \mu \cdot D^{k-1-i} \nabla \Phi_a ) \, dx   \end{bmatrix}_{a \in \{1,2,3\}} .
\end{eqnarray*}
Then we use that $\mathcal{M}$ is invertible and we apply the estimates  \eqref{pressurephi0},  \eqref{pressureclphi}, \eqref{pressuremudepart},  \eqref{pressureclmu} to obtain that there exists a positive decreasing function $\gamma_{4} $ with $\displaystyle {\lim_{L \rightarrow + \infty}} \gamma_{4}(L) =0$ such that 
\begin{equation}
\| \ell^{(k)} \| + | r^{(k)} | \leq   \gamma_{4}(L)   \sum_{i=0}^{k-1} \binom{k-1}{i}   \V_{i}  \V_{k-i-1} .
\end{equation}
Next we use  Lemma \ref{LemmeCheminSMF} in the case $s=2$ 
to get \eqref{MvtSolide}, where the function $ \gamma_{4}$ has been modified to incorporate the constant coming from \eqref{DefUpsilon}.
\ \par
In order to obtain \eqref{MvtFluide}, we write
\begin{equation}
D^{k} u = - D^{k-1} \nabla p = -D^{k-1} \nabla \mu + D^{k-1} \left(\nabla \Phi \cdot \begin{bmatrix}  \ell \\ r \end{bmatrix}'\right).	
	\label{new2.1}
\end{equation}
We notice that
$$
D^{k-1} \left(\nabla \Phi \cdot \begin{bmatrix}  \ell \\ r \end{bmatrix}'\right)
= \sum_{i=0}^{k-1} \binom{k-1}{i} D^i \nabla \Phi \cdot  \begin{bmatrix}  \ell^{(k-i)} \\ r^{(k-i)} \end{bmatrix}.
$$
Thus, by using \eqref{pressurehyp1} (valid up to rank $k$) and \eqref{pressureclphi} (valid up to rank $k$ due to Proposition \ref{pesti}) to estimate the terms of the above sum corresponding to $i\geq 1$ and by using \eqref{MvtSolide} for the term corresponding to $i=0$, we deduce, together with
\begin{gather*}
\sum_{i=1}^{k-1} \frac{1}{(i+1)^{2} (k-i+1)^{2}} \leq \sum_{i=1}^{k-1} \frac{1}{(i+1)^{2} }  \leq \frac{\pi^{2}}{6}, \\
\sum_{i=0}^{k-1} \binom{k-1}{i} \, (i!)^{M} ((k-i)!)^{M} \leq (k!)^{M},
\end{gather*}
that
$$
\left\|D^{k-1} \left(\nabla \Phi \cdot \begin{bmatrix}  \ell \\ r \end{bmatrix}'\right)\right\|\leq 
\left(\frac{\pi^{2}}{6} \gamma_1(L) + \gamma_4(L)C_0\right) \V_{k}.
$$
Combining the above inequality, \eqref{new2.1} and \eqref{pressureclmu}, we obtain \eqref{MvtFluide}, and the proof is complete.
\subsection{End of the proof }
\label{accord}
Let $\nu \in (0,1)$. We now apply  Eq. \eqref{proofit}  with  $p_2 (k+1)$ instead of $p_2$, and we use Stirling's formula to obtain that for any $k \in \N^*$, for any 
\begin{equation} \label{CondP2}
p_2  \geq \frac{2}{1-\nu},
\end{equation}
one has for $L$ sufficiently large (depending on the geometry only)
\begin{eqnarray}
\label{proofit3}
\| D^k u \|_{ W^{1,p_2 } (\mathcal{F}(t) ) }  + \| \ell^{(k)} \| + | r^{(k)} | 
& \leqslant  &  (p_2 (k+1))^k \frac{(k!)^{M} L^k }{(k+1)^2} 
\Big(  \|u\|_{ W^{1, p_2 (k+1)} (\mathcal{F}  (t) ) } +\|\ell\| + |r| \Big)^{k+1} \\
& \leqslant &  p_2^k (k!)^{M+1} \tilde{L}^k  \Big( \|u\|_{ W^{1, p_2 (k+1)} (\mathcal{F}  (t) ) } +\|\ell\| + |r| \Big)^{k+1}   ,
\end{eqnarray}
for some constant $\tilde{L} \geq L$ independent of $k$ and $(\ell,r,u)$. \par
So far time has intervened only as a parameter, and the inequality  \eqref{proofit3} holds for any time.
We will now estimate its right hand side with respect to the initial data.
First thanks to Lemma \ref{triv} there exists $c >0$ such that  for any $k$,
\begin{equation*}
\| u \|_{ W^{1,p_2  (k+1)} (\mathcal{F}(t) ) } \leq c p_2  (k+1)
\| \curl u \|_{ L^{p_2  (k+1)}(\mathcal{F}(t) ) } + c (| \gamma | + \| \ell \| + |r|).
\end{equation*}
Now conservation of the $L^p$ norms of the vorticity and Kelvin's circulation theorem yields that at any positive time
\begin{eqnarray*}
\| u \|_{ W^{1,p_2  (k+1)} (\mathcal{F}(t) ) } &\leq& c p_2  (k+1)
\| \omega_0 \|_{ L^{p_2  (k+1)}(\mathcal{F}(t) ) }  + c  (| \gamma | + \| \ell \| + |r|), \\
&\leq& c p_2  (k+1)   + c  (| \gamma | + \| \ell \| + |r|) ,
\end{eqnarray*}
since $\omega_0$ is bounded with compact support (enlarging $c$ if necessary).
Plugging this into  \eqref{proofit3} and using again Stirling's formula, we obtain that  there exists $L >0$ depending only on $\mathcal{F}(t) $ such that  for any $k$,
\begin{eqnarray*}
\| D^k u \|_{ W^{1,p_2 } (\mathcal{F}(t) )  }
\leqslant   p_2^k (k!)^{M+1} L^{k+1}   \big(k!  p_2^{k+1}  +  | \gamma  |^{k+1} +\| \ell \|^{k+1} + |r|^{k+1} \big) .
\end{eqnarray*}
Using  $p_2 \geq 2/(1-\nu)$, thanks to Morrey's inequality, there exists $C>0$ such that for any smooth function $f$ on $\overline{\mathcal{F}(t) }$, 
\begin{equation}
 \label{GraalUnif}
\| f \|_{ C^{0,\nu} (\mathcal{F}(t) )  } \leqslant  C  \| f \|_{ W^{1,p_2 } (\mathcal{F}(t))  }.
\end{equation}
This allows to bound $ \| D^k u \|_{ C^{0,\nu} (\mathcal{F}(t) ) }$ thanks to $\| D^k u \|_{ W^{1,p_2 } (\mathcal{F}(t)) }$. 
Then we differentiate Eq. \eqref{flow} to get 
\begin{equation*}
\partial^{k +1 }_t \Phi^\mathcal{F}  (t,x)= D^k u (t, \Phi^\mathcal{F} (t,x)).
\end{equation*}
We consider $T>0$ and we obtain by composition that $\partial^{k+1}_t \Phi^\mathcal{F}$ is in $C^{0, \nu \exp(-c T \| \omega_0 \|_{L^{\infty}  ({\mathcal F}_0)})} ({\mathcal F}_0)$, with an estimate
\begin{equation} \label{Graal}
\| \partial^{k+1}_t \Phi^\mathcal{F} \|_{ C^{0, \nu \exp(-c T \| \omega_0 \|_{L^{\infty}  ({\mathcal F}_0)})} ({\mathcal F}_0) }
\leqslant   p_2^k (k!)^{M+1} L^{k+1}   \big(k!  p_2^{k+1}  +  | \gamma  |^{k+1} +\| \ell \|^{k+1} + |r|^{k+1} \big) ,
\end{equation}
for some $L$ depending on the geometry only. \par
Now, in order to prove \eqref{PPPM}, we have to absorb the $\nu$ factor in the H\"older exponent in \eqref{Graal}.  To do this we use the fact that these estimates are valid whatever the choice of the time interval and of $\nu \in (0,1)$ (note that $p_{2}$ satisfies \eqref{CondP2} and that the constant in \eqref{GraalUnif} is uniform for $\nu \in (1/2,1)$). In particular, we consider $\tau \in (0,T)$, and we apply the above result on the interval $(-\tau , \tau)$ and $\nu \in (0,1)$ such that $\nu  \exp(-c  \tau \| \omega_0 \|_{L^{\infty} ({\mathcal F}_0)})  > \exp(-c T \| \omega_0 \|_{L^{\infty} ({\mathcal F}_0)})$. Then one can choose the constant $L>0$ containing the $(\cdot)^{k}$ factors and the proof of Theorem \ref{start4} is over. \par
\ \par
{\bf Acknowledgements.} The  authors were partially supported by the Agence Nationale de la Recherche, Project CISIFS,  grant ANR-09-BLAN-0213-02. 

\end{document}